\documentclass[final,3p,times,authoryear]{elsarticle}
\usepackage{hyperref}
\hypersetup{
    bookmarks=true,         
    unicode=false,          
    pdftoolbar=true,        
    pdfmenubar=true,        
    pdffitwindow=false,     
    pdfstartview={FitH},    
    pdftitle={The Algorithm of Factor Twelve},    
    pdfauthor={F. Quinonez},     
    pdfsubject={MSC: 01A17},   
    pdfcreator={F. Quinonez},   
    pdfproducer={F. Quinonez}, 
    pdfkeywords={Plimpton 322, Old Babylonian Trigonometry}, 
    pdfnewwindow=true,      
    colorlinks=true,       
    linkcolor=magenta,          
    citecolor=green,        
    filecolor=red,      
    urlcolor=cyan           
}

\usepackage{pdflscape}
\usepackage{eqnarray}
\usepackage{amsmath, amsthm, amssymb, amsfonts}
\usepackage{graphicx}

\usepackage{mathscinet}
\ProvideTextCommandDefault{\rasp}{\leavevmode\raise.45ex\hbox{$\rhook$}}
\ProvideTextCommandDefault{\lasp}{\leavevmode\raise.45ex\hbox{$\lhook$}}

\newcommand{\divides}{\mid}

\newcommand {\zero}{\ensuremath{\:\vcenter{\hbox{\includegraphics[height=5mm]{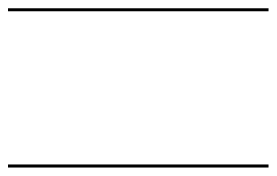}}}}}
\newcommand {\one}{\ensuremath{\:\:\vcenter{\hbox{\includegraphics[height=5mm]{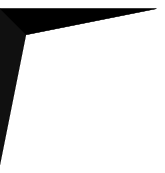}}}}}
\newcommand {\two}{\ensuremath{\:\:\vcenter{\hbox{\includegraphics[height=5mm]{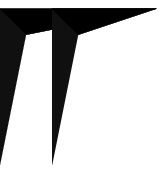}}}}}
\newcommand {\three}{\ensuremath{\:\:\vcenter{\hbox{\includegraphics[height=5mm]{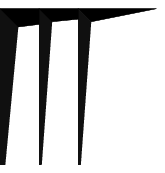}}}}}
\newcommand {\four}{\ensuremath{\:\:\vcenter{\hbox{\includegraphics[height=5mm]{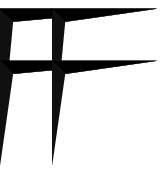}}}}}
\newcommand {\five}{\ensuremath{\:\:\vcenter{\hbox{\includegraphics[height=5mm]{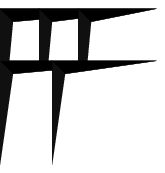}}}}}
\newcommand {\six}{\ensuremath{\:\:\vcenter{\hbox{\includegraphics[height=5mm]{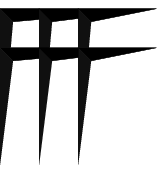}}}}}
\newcommand {\seven}{\ensuremath{\:\:\vcenter{\hbox{\includegraphics[height=5mm]{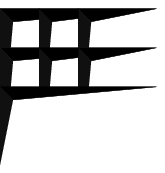}}}}}
\newcommand {\eight}{\ensuremath{\:\:\vcenter{\hbox{\includegraphics[height=5mm]{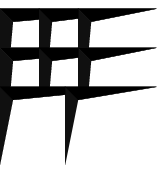}}}}}
\newcommand {\nine}{\ensuremath{\:\:\vcenter{\hbox{\includegraphics[height=5mm]{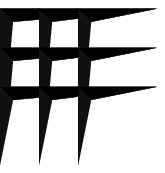}}}}}
\newcommand {\ten}{\ensuremath{\vcenter{\hbox{\includegraphics[height=5mm]{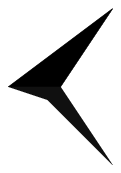}}}}}
\newcommand {\eleven}{\ensuremath{\vcenter{\hbox{\includegraphics[width=9mm]{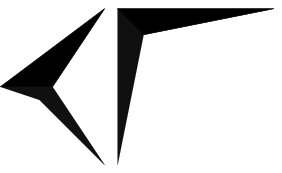}}}}}
\newcommand {\twelve}{\ensuremath{\vcenter{\hbox{\includegraphics[width=9mm]{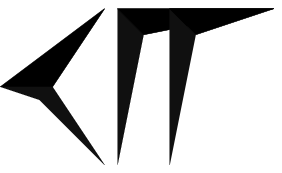}}}}}
\newcommand {\thirteen}{\ensuremath{\vcenter{\hbox{\includegraphics[width=9mm]{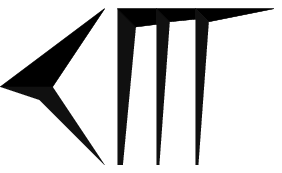}}}}}
\newcommand {\fourteen}{\ensuremath{\vcenter{\hbox{\includegraphics[width=9mm]{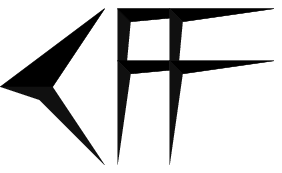}}}}}
\newcommand {\fifteen}{\ensuremath{\vcenter{\hbox{\includegraphics[width=9mm]{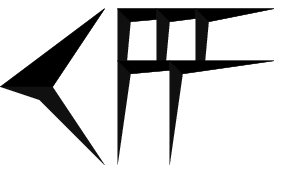}}}}}
\newcommand {\sixteen}{\ensuremath{\vcenter{\hbox{\includegraphics[width=9mm]{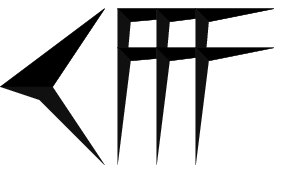}}}}}
\newcommand {\nineteen}{\ensuremath{\vcenter{\hbox{\includegraphics[width=9mm]{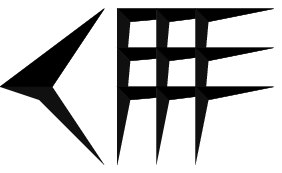}}}}}
\newcommand {\twenty}{\ensuremath{\vcenter{\hbox{\includegraphics[height=5mm]{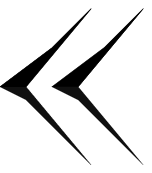}}}}}
\newcommand {\twentyone}{\ensuremath{\vcenter{\hbox{\includegraphics[width=10mm]{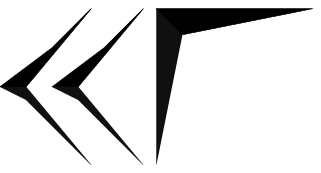}}}}}
\newcommand {\twentytwo}{\ensuremath{\vcenter{\hbox{\includegraphics[width=10mm]{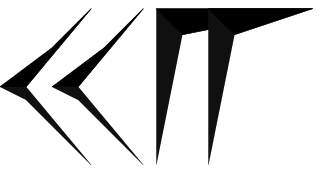}}}}}
\newcommand {\twentythree}{\ensuremath{\vcenter{\hbox{\includegraphics[width=10mm]{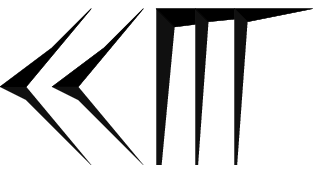}}}}}
\newcommand {\twentyfour}{\ensuremath{\vcenter{\hbox{\includegraphics[width=10mm]{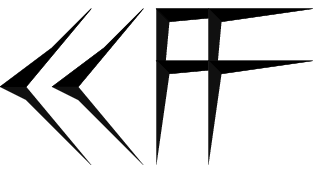}}}}}
\newcommand {\twentyfive}{\ensuremath{\vcenter{\hbox{\includegraphics[width=10mm]{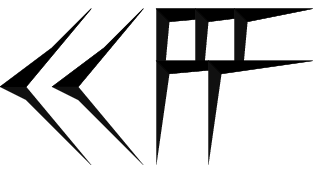}}}}}
\newcommand {\twentysix}{\ensuremath{\vcenter{\hbox{\includegraphics[width=10mm]{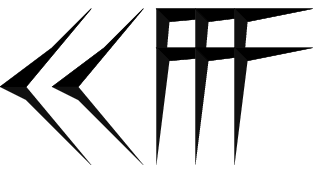}}}}}
\newcommand {\twentyseven}{\ensuremath{\vcenter{\hbox{\includegraphics[width=10mm]{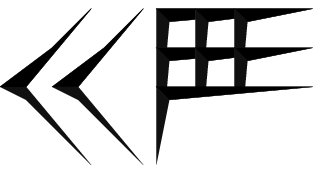}}}}}
\newcommand {\twentyeight}{\ensuremath{\vcenter{\hbox{\includegraphics[width=10mm]{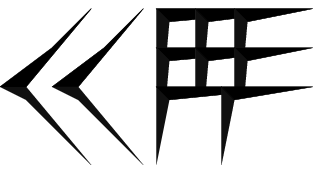}}}}}
\newcommand {\twentynine}{\ensuremath{\vcenter{\hbox{\includegraphics[width=10mm]{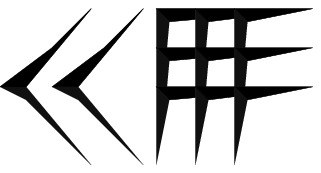}}}}}
\newcommand {\thirtyone}{\ensuremath{\vcenter{\hbox{\includegraphics[width=10mm]{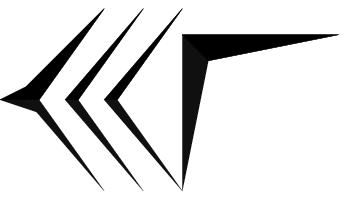}}}}}
\newcommand {\thirtytwo}{\ensuremath{\vcenter{\hbox{\includegraphics[width=10mm]{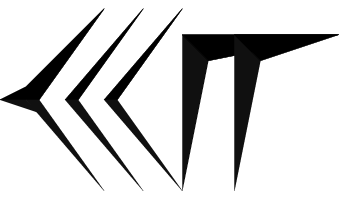}}}}}
\newcommand {\thirtythree}{\ensuremath{\vcenter{\hbox{\includegraphics[width=10mm]{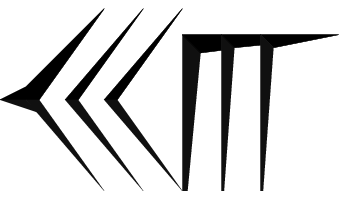}}}}}
\newcommand {\thirtyfive}{\ensuremath{\vcenter{\hbox{\includegraphics[width=10mm]{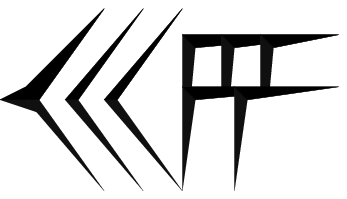}}}}}
\newcommand {\thirtysix}{\ensuremath{\vcenter{\hbox{\includegraphics[width=10mm]{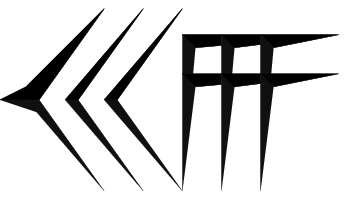}}}}}
\newcommand {\thirtyseven}{\ensuremath{\vcenter{\hbox{\includegraphics[width=10mm]{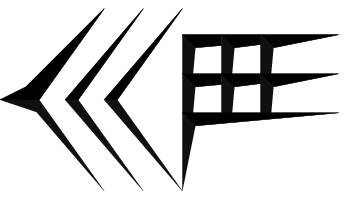}}}}}
\newcommand {\thirtyeight}{\ensuremath{\vcenter{\hbox{\includegraphics[width=10mm]{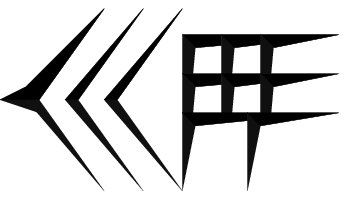}}}}}
\newcommand {\fourty}{\ensuremath{\vcenter{\hbox{\includegraphics[width=5mm]{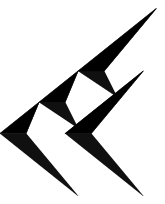}}}}}
\newcommand {\fourtyone}{\ensuremath{\vcenter{\hbox{\includegraphics[width=10mm]{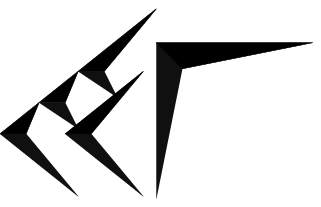}}}}}
\newcommand {\fourtythree}{\ensuremath{\vcenter{\hbox{\includegraphics[width=10mm]{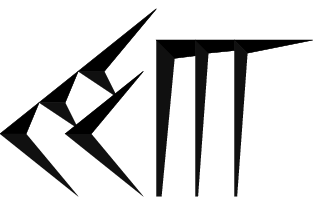}}}}}
\newcommand {\fourtyfive}{\ensuremath{\vcenter{\hbox{\includegraphics[width=10mm]{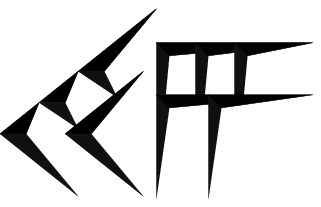}}}}}
\newcommand {\fourtysix}{\ensuremath{\vcenter{\hbox{\includegraphics[width=10mm]{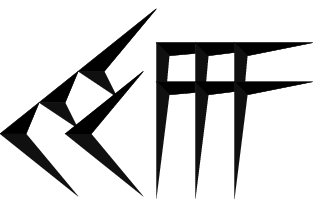}}}}}
\newcommand {\fourtyseven}{\ensuremath{\vcenter{\hbox{\includegraphics[width=10mm]{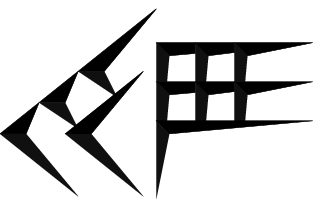}}}}}
\newcommand {\fourtyeight}{\ensuremath{\vcenter{\hbox{\includegraphics[width=10mm]{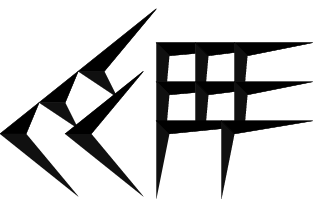}}}}}
\newcommand {\fourtynine}{\ensuremath{\vcenter{\hbox{\includegraphics[width=10mm]{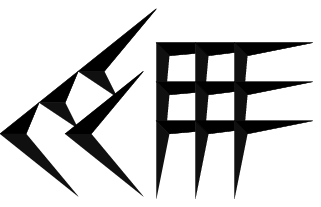}}}}}
\newcommand {\fifty}{\ensuremath{\vcenter{\hbox{\includegraphics[width=5mm]{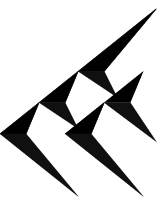}}}}}
\newcommand {\fiftyone}{\ensuremath{\vcenter{\hbox{\includegraphics[width=10mm]{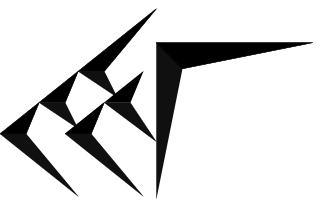}}}}}
\newcommand {\fiftytwo}{\ensuremath{\vcenter{\hbox{\includegraphics[width=10mm]{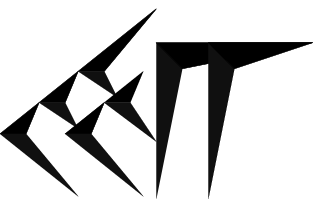}}}}}
\newcommand {\fiftythree}{\ensuremath{\vcenter{\hbox{\includegraphics[width=10mm]{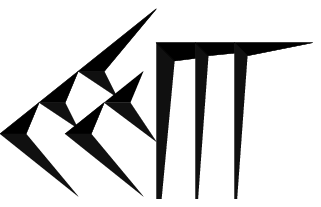}}}}}
\newcommand {\fiftyfour}{\ensuremath{\vcenter{\hbox{\includegraphics[width=10mm]{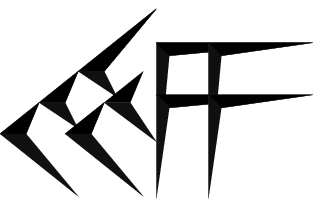}}}}}
\newcommand {\fiftyfive}{\ensuremath{\vcenter{\hbox{\includegraphics[width=10mm]{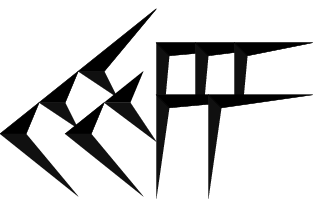}}}}}
\newcommand {\fiftysix}{\ensuremath{\vcenter{\hbox{\includegraphics[width=10mm]{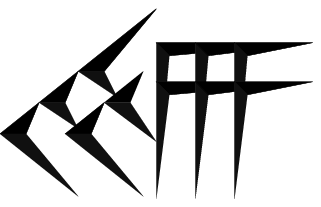}}}}}
\newcommand {\fiftyeight}{\ensuremath{\vcenter{\hbox{\includegraphics[width=10mm]{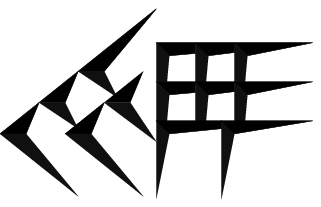}}}}}
\newcommand {\fiftynine}{\ensuremath{\vcenter{\hbox{\includegraphics[width=10mm]{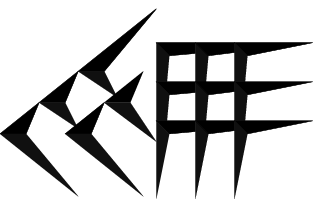}}}}}

\newcommand {\cfive}{\ensuremath{[05]}}
\newcommand {\cfifteen}{\ensuremath{[15]}}
\newcommand {\ctwentyfive}{\ensuremath{[25]}}
\newcommand {\cfourtyfive}{\ensuremath{[45]}}
\newcommand {\cseven}{\ensuremath{[07]}}
\newcommand {\celeven}{\ensuremath{[11]}}
\newcommand {\ctwelve}{\ensuremath{[12]}}
\newcommand {\cthirteen}{\ensuremath{[13]}}
\newcommand {\cseventeen}{\ensuremath{[17]}}
\newcommand {\cnineteen}{\ensuremath{[19]}}
\newcommand {\ctwentythree}{\ensuremath{[23]}}
\newcommand {\ctwentyeight}{\ensuremath{[28]}}
\newcommand {\ctwentynine}{\ensuremath{[29]}}
\newcommand {\cthirtyone}{\ensuremath{[31]}}
\newcommand {\cthirtysix}{\ensuremath{[36]}}
\newcommand {\cthirtyseven}{\ensuremath{[37]}}
\newcommand {\cfourtyone}{\ensuremath{[41]}}
\newcommand {\cfourtythree}{\ensuremath{[43]}}

\newcommand {\cfourtynine}{\ensuremath{[49]}}
\newcommand {\cfiftythree}{\ensuremath{[53]}}
\newcommand {\cfiftynine}{\ensuremath{[59]}}
\newcommand {\cone}{\ensuremath{[01]}}
\newcommand {\czero}{\ensuremath{[00]}}
\newcommand {\cfourtyeight}{\ensuremath{[49]}}

\begin{document}
\begin{frontmatter}

\title{The Algorithm of Factor 12: \\Generating the Information Carved on the Plimpton 322 Tablet.}
%
\author{Fernando Qui\~{n}onez}
\ead{fquinonezf@gmail.com, faquinon@uc.cl}

%
%
%
%
\begin{abstract}
The right triangles represented in the Plimpton 322 data have integer sides $(a, b, d)$ with
$a < b < d$, by inspecting the data we propose $b=M\,Q_{M}$ with $M,Q_{M}$ positive integers. 
In this way we present a new data-driven fa\-mi\-ly of algorithms in order to generate
the whole content of the Plimpton-322 tablet.
Each algorithm of the fa\-mi\-ly co\-rres\-ponds to one given va\-lu\-e of the bundling factor or \emph{mak\c{s}arum} $M$.
The algorithm with $M=12$ arise naturally from data 
and is developed entirely through this work. 
We could generate the complete Plimpton 322 table by using three different schemes 
in the framework of the twelve factor algorithm.
Finally we found by varying $M$ in the general algorithm, we can be able to find all the Pythagorean Triples.
\end{abstract}
\end{frontmatter}

\section{Introduction}
The 322nd clay tablet in the Plimpton collection of Columbia University in New York has
wondered historian-mathematicians since \citep{Neugebauer1935} 
and also \citep{NeugebauerSachsGoetze1945} 
revealed this tablet contains carved information about fifteen right triangles. 
The Plimpton-322 tablet (P322) is an ancient artifact with rectangular shape of height $8.8\;\mathrm{cm}$ and width $12.7 \;\mathrm{cm}$,
each row has a height of $5 \;\mathrm{mm}$ or $6\;\mathrm{mm}$ and are delimited by guide lines carved horizontally
through the obverse of the entire tablet.
As a matter of fact, there are not separation between the lines that delimit each row 
and some sexagesimal cuneiform numbers go out of its row field mostly invading the next row that comes above it.
This forced the scribes to left blank spaces between some numbers trying to prevent some invasion from the next row below it and 
like in all preventive actions, sometimes they works and sometimes not, that is the reason 
that not all the invasions were supported by blank spaces and also not all blank spaces between numbers in the same column were invaded.
P322 triangles are plotted scaled in Figure \ref{f:rt} and 
we have a representation of P322 in Table \ref{t:Plimpton322}
but we miss the guide lines and we have the default row separation given by the \LaTeX\ tables.
The parenthesis enclose pieces of the tablet that are damaged whilst the square parenthesis enclose the errors made by the scribe.
To look the real P322 see \citep{RobsonAMSoutreach}. P322 is dated by \citep{Robson1999}, \citep{Robson2001}, to belong 
to the period between 1822 and 1762 B.C.
(3842-3782 years ago) and therefore belongs to the Old-Babylonian period ($[1800, 1600]$ B.C.). 
This period is mainly characterized by the non presence of a symbol to represent the number zero,
without means they did not know the number zero, 
instead they used a blank space too, adding more chaos to our analysis. 
However the combined efforts of many authors along the years have helped us to understand how the P322 information was filled and at the end, 
we can find in the literature there are three mainstream theories: 
the early proposal by Neugerbauer based on the Euclid transformation 
for Pythagorean triangles triples of the form
$(p^{2} - q^{2}, 2pq, p^{2} + q^{2})$ and then he was able to find the left column as $(p^{2} + q^{2})^{2} / (2pq)^{2}$,
the proposed reciprocal theory \citep{Bruins1949} $((x - 1/x)/2, 1, (x + 1/x)/2)$ with left column 
$(x - 1/x)^{2}/4$
and adopted by \citep{Robson1999} with left column
$(x + 1/x)^{2}/4$, 
and the third is a family of little modifications of the former Neugebauer proposal 
first made by 
\citep{DeSollaPrice1964}, 
\citep{Aaboe1964},
$((p\bar{q} - q\bar{p})/2, 1, (p\bar{q} + q\bar{p})/2)$
and followed by the works of \citep{Friberg1981,Friberg2000,Friberg2001}, 
\citep{Proust2002}, 
\citep{Hoyrup2002}, %
\citep{BrittonProustSchnider2011},
\citep{Abdulaziz2010},
and finally \citep{Mansfield2017}, they arrive to something like
$((\frac{1}{2}(r\bar{s} - s\bar{r}), 1, \frac{1}{2}(r\bar{s} + s\bar{r})))$.
All can be considered as generating algorithms able to reconstruct the nu\-me\-ri\-cal 
values presented in P322. 
\\\\\indent Their works are supported on a set of 
techniques that are attributed to the scribes of the time, some of these have been 
proved by the evidence found in other tablets with cuneiform carved on them.
P322 was found in Larsa in a zone where another numerical tablets known as Multiplication Tables (MT) have been found too. 
Many others stone tablets have been lost by the pass of time and another ones have been analyzed and found not much useful as a simple (MT),
however there is one MT with a more deep meaning, it is: 
$\one \twelve = 01.~12 = 6/5$.
Examples of MT can be seen in
\citep{AMS-MT}: 
VAT7858 a ten-times table, 
NBC7344 a five-times table,
MS2184/3 and NBC7344 as twelve-times tables,
MS3866 1.2 or 72 times table, 
and 
MS3874 a table of reciprocals.
%
\subsection{The Amazing Old Babylonian Place Value Notation.}
As Neugebauer explains in his classic book \citep{Neugebauer1957}, P322 presents a technique called
the Old Babylonian place value notation. 
Basically it means as the Old Babylonian scribes does not use a point or another symbol to represent
where starts the floating part of a number, or said in modern mathematical terms, they
could interpret freely where to put the sexagesimal point.
As this notation implicitly involves multiplication of the number times the base raised to an integer exponent $60^{n}$.
In scientific terms, one could interpret this notation to be one of the earliest implicit use of the \emph{scientific notation}. 
\\\\\indent Right triangles triples and the left column in P322 share another property, suppose the triple is
$(a,b,d)$ it satisfies $a^{2}+b^{2}=d^{2}$ and the left column is $d^{2}/b^{2}$ or $a^{2}/b^{2}$, now if one has
$(ka,kb,kd)$ with $k$ integer and not necessarily of the form $60^{n}$, the new scaled triple and the left column will satisfy the same condition and value. 
Yes it is a geometric triviality but is a key property used in this work.
\\\\\indent This work presents a new algorithm to generate the data of P322,
based on the P322 data itself, and the
set of MT tables, square/square-root tables and the magnificent Reciprocals Table (RT) 
found near in that place-age,
together with the Old Babylonian place value notation;
this is a data driven algorithm in the context of that age.
\section{Data Description}
\begin{figure}[h!]
\centering
\includegraphics[width=12cm]{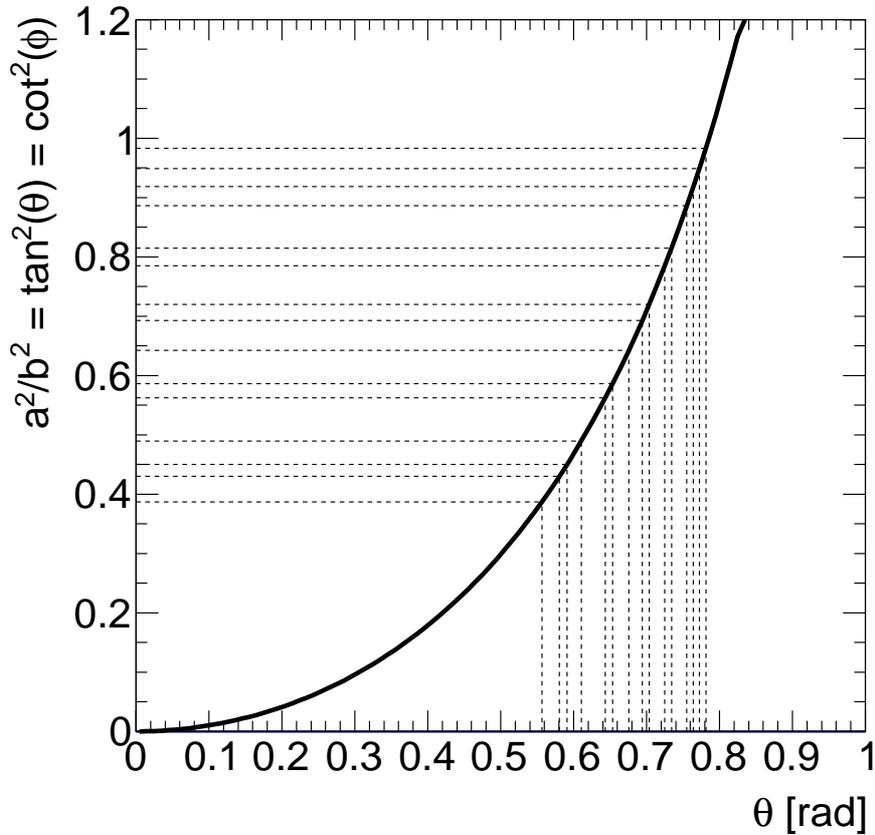}
\caption{Values for the most to the left column in P322, they can be either $a^{2}/b^{2}$ or $d^{2}/b^{2} = 1 + a^{2}/b^{2} = \sec^{2}{\theta}$.}
\label{f:leftcolumn}
\end{figure}

In Figure~\ref{t:Plimpton322} the
true labels (titles) of columns are omitted and we have replaced them by labels,
from right to left: 
`its name', $d$, $a$, and $d^{2}/b^{2}$ or $a^{2}/b^{2}$; instead of the original Babylonian symbols that according to 
\citep[p. 107]{Robson2002} they are
MU.BI.IM, \'IB.SI$_{8}$ {\it \d{s}i-li-ip-tim} ({\it mit\uarc{h}arti \d{s}iliiptim}), \'IB.SI$_{8}$ SAG ({\it mit\uarc{h}arti p\={u}tim}), and 
{\it ta-ki-il-ti \d{s}i-li-ip-tim \v{s}a} 1 {\it in-na-as-s\`{a}-\uarc{h}u-ma} SAG {\it i-il-lu-\'{u}};
that means:
as itself, the square side of the diagonal, the square side of the short side, and  
the title of the 4th column is according to 
\citep[p. 191]{Robson2002} and the references there: 
\emph{``The takiltum of the diagonal from which 1 is torn out, so that the short side comes up''};
respectively.
\begin{figure}
\centering
\includegraphics[width=7.6cm]{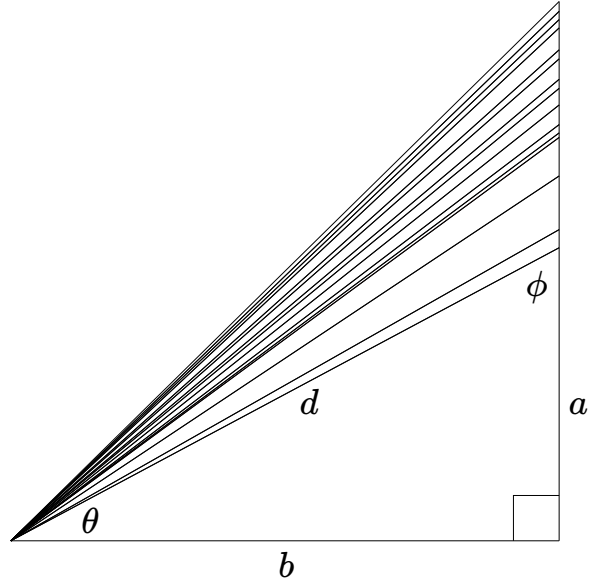}
\caption{Scaled right triangles representation of the fifteen right triangles in P322. They are ordered from top to bottom according to the same order in P322.}
\label{f:rt}
\end{figure}
\\\\\indent Every triangle in P322 has different integer sides $a < b < d$.
\\\\\indent Every entry row in the 1st column is a sexagesimal number with va\-lu\-es:
one ($\one = 01$), 
two ($\two = 02$),  
three ($\three = 03$),  
four ($\four = 04$),  
five ($\five = 05$),  
six ($\six = 06$),  
seven ($\seven = 07$),  
eight ($\eight = 08$),  
nine ($\nine = 09$),  
ten ($\ten = 10$),  
eleven ($\eleven = 11$),  
twelve ($\twelve = 12$),  
thirteen ($\thirteen = 13$),  
fourteen ($\fourteen = 14$),  
and fifteen ($\fifteen = 15$). 
Before each number there is a word \emph{ki} 
that according to the experts means \emph{``as itself'}.
\\\\\indent A very remarkable feature of P3222 is that the numbers in the 4th column comes sorted from highest to lowest.
\\\\\indent In 4th column the value $a^{2}/b^{2}$ or equivalently $d^{2}/b^{2}-1$ can be seen in the ordinate axis for $\tan^{2}{\theta}$ 
in Figure~\ref{f:leftcolumn}, where
$\theta$ is the shortest angle between the sides $b$ and $d$, see $\theta$ in Figure~\ref{f:rt}. 
\begin{table}[t!]
\centering
\begin{tabular}{llll}
\hline
$\frac{d^{2}}{b^{2}} \big| \frac{d^{2}}{b^{2}} - 1 = \frac{a^{2}}{b^{2}}$ & $a$ & $d$ & mu bi im \\
\hline
column 4 & column 3 & column 2 & column 1 \\
$(\one\fiftynine) \zero\fifteen$  &  \one\fiftynine & \two \fourtynine & ki \one \\
$(\one\fiftysix\fiftysix) \fiftyeight\fourteen\fiftysix\fifteen$  & \fiftysix\seven & [\three\twelve\one ] & ki \two \\
$(\one\fiftyfive\seven) \fourtyone\fifteen\thirtythree\fourtyfive$   & \one\sixteen\fourtyone    &\one\fifty\fourtynine  & ki \three \\
$(\one\fiftythree\ten) \twentynine\thirtytwo\fiftytwo\sixteen$  & \three\thirtyone\fourtynine & \five\nine\one & ki \four \\
$(\one) \fiftyeight\fiftyfour\zero\one\fourty$   &  \one\five   & \one\thirtyseven   & ki (\five) \\
$(\one) \fourtyseven\zero\six\fourtyone\fourty$   &   \five\nineteen    &  \eight\one & (ki \six) \\
$(\one) \fourtythree\eleven\fiftysix\twentyeight\twentysix \fourty$   &  \thirtyeight\eleven & \fiftynine\one & ki \seven \\
$(\one) \fourtyone\thirtythree\fiftynine    \zero   \three\fourtyfive$   &  \thirteen \nineteen    &  \twenty\fourtynine  & ki \eight  \\
$(\one) \thirtyeight\thirtythree\thirtysix\thirtysix$   &  [\eight ]\one   &  \twelve\fourtynine   &  ki \nine  \\
$(\one) \thirtyfive\ten\two\twentyeight\twentyseven\twentyfour\twentysix\fourty$  & \one\twentytwo\fourtyone &\two\sixteen\one & ki \ten \\
$(\one) \thirtythree\fourtyfive$   &  \fourtyfive\zero    & \one\fifteen\zero  &  ki \eleven  \\
$(\one) \twentynine\twentyone\fiftyfour   \zero   \two\fifteen$   &  \twentyseven\fiftynine & \fourtyeight\fourtynine& ki  \twelve  \\
$(\one) \twentyseven   \zero   \three\fourtyfive$   & [\seven\twelve\one ]    &  \four\fourtynine &  ki \thirteen  \\
$(\one) \twentyfive\fourtyeight\fiftyone\thirtyfive \zero  \six\fourty $& \twentynine\thirtyone & \fiftythree\fourtynine & ki \fourteen  \\
$(\one) \twentythree\thirteen\fourtysix\fourty$   &   \fiftysix     &  \fiftythree   &  ki (\fifteen)  \\
\hline
\end{tabular}
\caption{P322 Representation. This picture keeps the scribe's errors which are enclosed between square parenthesis, 
normal parenthesis means these parts of the tablet are damaged.}
\label{t:Plimpton322}
\end{table}
\\\indent Another remarkable and controversial feature in these numbers comes  
from the point of view of the $\theta$ angle, this is, it is bounded.
Said in modern terms: $\theta \in (\pi/4, \pi/6)$, avoiding the extremes of the interval. 
If the intervals were taken, i.e. $\theta=\pi/4$ will give right triangles with triples 
$( a, a, a \sqrt{2} )$ 
and $\theta=\pi/6$ will give $( a, a \sqrt{3}, 2a)$, implying 
the use of infinite sexagesits in each case in order to represent the exact value of 
$\sqrt{2} \cdot a$ and $\sqrt{3} \cdot a$.
Some authors do not consider the P322 as a trigonometrical table and arguing the Old Babylonians did not know $\pi$ and that 
there was not explicit the scribe was varying and angle.
\\\\\indent Let $s =_{10} d$ denote the sexagesimal number $s$ represented in decimal base by the number $d$, then 
the detailed content by rows is as follows:
\\\\{\it Row 1:} $d=\two \fourtynine = 02~49 =_{10} 169=13^{2}$, $a=\one\fiftynine=01~59 =_{10} 119 = 7\cdot17$, 
$a^{2}/b^{2} = \fiftynine\zero\fifteen = 59~00~15 =_{10} 212415 \times 60^{-3}$ or
$d^{2}/b^{2} = \one\fiftynine\zero\fifteen = 01~59~00~15 =_{10} 428415 \times 60^{-3}$.
\\\\{\it Row 2:} This row has an error made by the scribe. The thing is the 4th column has at least two sexagesits damaged 
and the entire number can not be readed. 
But if one maintain fixed the last sexagesits and vary the two
unknown we note easily the triple presents an error, because we can be able to find two sexagesits for the number $\frac{a^{2}}{d^{2}-a^{2}}$. 
Also if one maintain fixed $d$ but varying $a$ we can not be able to find the solution again.
However, as noted by \citep[p. 38]{Neugebauer1957} if one change the value of the diagonal 
from $d=\three\twelve\one = 03~12~01 =_{10} 11521 = 41 \cdot 281$
to be 
$d=\one \twenty \twentyfive = 01~20~25 =_{10} 4825 = 5^{2} \cdot 193$, the triple is fixed and we have $a=\fiftysix\seven= 56~07 =_{10} 3367 = 7 \cdot 13 \cdot 37$, 
$a^{2}/b^{2} = \fiftysix \fiftysix \fiftyeight \fourteen \fiftysix \fifteen = 56~56~58~14~56~15 =_{10} 44283941775\times60^{-6}$ or
$d^{2}/b^{2} = 01~56~56~58~14~56~15 =_{10} 90939941775\times60^{-6}$.
About this error Neugebauer said: 
\emph{``Finally there remains unexplained error in III,2 where 3,12,1 should be replaced by 1,20,25''}. In \citep[p. 193]{Robson2001} explains this error
not exactly.
For instance in this work we explain the error by noting
$$ 03~12~01 = (02~41)^{2} - (02~00)^{2} $$
is an error related by operating with $b=_{10} 120$ coming from Row 1:
$(01~59)^{2} = (02~49)^{2} - (02~00)^{2}$.
\\\\{\it Row 3:} $d=\one \fifty \fourtynine = 01~50~49 =_{10} 6649 = 61\cdot109$, $a=\one \sixteen \fourtyone = 01~16~41 =_{10} 4601 = 43\cdot107$, 
$a^{2}/b^{2} = \fiftyfive\seven\fourtyone \fifteen \thirtythree \fourtyfive = 55~07~41~15~33~45 =_{10} = 428676632025\times 60^{-6}$ or
$d^{2}/b^{2} = 01~55~07~41~15~33~45 =_{10} = 89523632025\times 60^{-6}$.
\\\\{\it Row 4:} $d=\five\nine\one = 05~09~01  =_{10} 18541$ (prime), $a=\three\thirtyone\fourtynine = 03~31~49 =_{10} 12709= 71\cdot179$, 
$a^{2}/b^{2} = \fiftythree\ten\twentynine\thirtytwo\fiftytwo\sixteen = 53~10~29~32~52~16 =_{10} 41348782336\times60^{-6}$ or
$d^{2}/b^{2} = 01~53~10~29~32~52~16 =_{10} 88004782336\times60^{-6}$.
\\\\{\it Row 5:} $d=\one \thirtyseven = 01~37  =_{10} 97$ (prime), $a=\one\five = 01~05 =_{10} 72$, 
$a^{2}/b^{2} = \fiftyeight\fiftyfour\one\fourty = 58~54~01~40 =_{10} 12722500\times60^{-4}$ or
$d^{2}/b^{2} = 01~58~54~01~40 =_{10} 25682500\times60^{-4}$.
\\\\{\it Row 6:} $d=\eight\one = 08~01  =_{10} 481 = 13\cdot37$, $a=\five\nineteen = 05~19 =_{10} 319 = 11\cdot29$, 
$a^{2}/b^{2} = \fourtyseven \six \fourtyone \fourty = 47~06~41~40 =_{10} 10176100\times60^{-4}$ or
$d^{2}/b^{2} = 01~47~06~41~40 =_{10} 23136100\times60^{-5}$.
\\\\{\it Row 7:} $d=\fiftynine\one = 59~01  =_{10} 3541$ (prime), $a=\thirtyeight\eleven = 38~11 =_{10} 2291 = 29\cdot79$, 
$a^{2}/b^{2} = \fourtythree\eleven\fiftysix\twentyeight\twentysix\fourty = 43~11~56~28~26~40 =_{10} 33591558400\times60^{-6}$ or
$d^{2}/b^{2} = 01~43~11~56~28~26~40 =_{10} 80427558400\times60^{-6}$.
\\\\{\it Row 8:} $d=\twenty\fourtynine = 20~49  =_{10} 1249$ (prime), $a=\thirteen\nineteen = 13~19 =_{10} 799 = 17\cdot47$, 
$a^{2}/b^{2} = \fourtyone \thirtythree\fourtyfive\fourteen\three\fourtyfive = 41~33~45~14~03~45 =_{10} 32319050625\times60^{-6}$ or
$d^{2}/b^{2} = 01~41~33~45~14~03~45 =_{10} 78975050625\times60^{-6}$.
\\\\{\it Row 9:} In this row there is an error typo, also noted by Neugebauer, in the value of the short side from (prime) $a=\nine\one = 09~01 =_{10} 541$ 
to $a=\eight\one= 08~01 =_{10} 481 = 13\cdot 37$, then 
$d=\twelve\fourtynine = 12~49  =_{10} 769$ (prime),  
$a^{2}/b^{2} = \thirtyeight\thirtythree\thirtysix\thirtysix = 38~33~36~36 =_{10} 8328996\times60^{-4}$ or
$d^{2}/b^{2} = 01~38~33~36~36 =_{10} 21288996\times60^{-4}$.
\\\\{\it Row 10:} $d=\two\sixteen\one = 02~16~01  =_{10} 8161$ (prime), $a=\one\twentytwo\fourtyone = 01~22~41 =_{10} 4961 = 11^{2}\cdot 41$, 
$a^{2}/b^{2} =$ $\thirtyfive\ten\two\twentyeight\twentyseven\twentyfour\twentysix\fourty=35~10~02~28~27~24~26~40=_{10}98446084000000\times60^{-8}$ or
$d^{2}/b^{2} = 01~35~10~02~28~27~24~26~40=_{10} 266407684000000\times60^{-8}$.
\\\\{\it Row 11:} $d=\one\fifteen\zero = 01~15~00  =_{10} 4500 = 2^{2}\cdot 3^{2} \cdot 5^{3}$, $a=\fourtyfive\zero =_{10} 2700= 2^{2}\cdot 3^{3} \cdot 5^{2}$, 
$a^{2}/b^{2} = \thirtythree\fourtyfive = 33~45 =_{10} 2025\times60^{-2}$ or
$d^{2}/b^{2} = 01~33~45 =_{10} 5625\times60^{-2}$.
\\\\{\it Row 12:} $d=\fourtyeight\fourtynine = 48~49  =_{10} 2929 = 29\cdot101$, $a=\twentyseven\fiftynine = 27~59 =_{10} 1679 = 23\cdot73$, 
$a^{2}/b^{2} = \twentynine\twentyone\fiftyfour\two\fifteen = 29~21~54~02~15 =_{10} 380570535\times60^{-5}$ or
$d^{2}/b^{2} = 01~29~21~54~02~15 =_{10} 1158170535\times60^{-5}$.
\\\\{\it Row 13:} Here there is another error and is due to the scribe squared the true value of the short side, 
again Neugebauer noted
$$ 07~12~01 = (02~41)^{2} $$
then we have to replace $a=\seven\twelve\one = 07~12~01 =_{10} 25921$ to its square-root 
$a=\two\fourtyone =_{10} 161= 7\cdot23$, then
$d=\four\fourtynine = 04~49  =_{10} 289 = 17^{2}$, 
$a^{2}/b^{2} = \twentyseven\zero\three\fourtyfive = 27~00~03~45 =_{10} 5832225\times60^{-4}$ or
$d^{2}/b^{2} = 01~27~00~03~45 =_{10} 18792225\times60^{-4}$.
\\\\{\it Row 14:} $d= \fiftythree\fourtynine = 53~49  =_{10} 3229$ (prime), $a=\twentynine\thirtyone = 29~31 =_{10} 1771 = 7\cdot 11 \cdot 23$, 
$a^{2}/b^{2} =\twentyfive\fourtyeight\fiftyone\thirtyfive\six\fourty = 25~48~51~35~06~40 =_{10} 20073222400\times60^{-6}$.
$d^{2}/b^{2} = 01~25~48~51~35~06~40 =_{10} 66729222400\times60^{-6}$.
\\\\{\it Row 15:} The original values $d=\fiftythree = 53  =_{10} 53$ (prime), $a=\fiftysix =_{10} 56$, 
fixing this could be done in several ways:
\\\indent{\it a la} Neugebauer: doubling the diagonal to get $d=01~46$,
\\\indent{\it a la} Robson: halfing the short side to get $a=28$,
\\\indent or maybe take $d=53~00$, $a=28~00$; 
\\but we will take the third option at the moment;
then $a^{2}/b^{2} = \twentythree\thirteen\fourtysix\fourty  = 23~13~46~40 =_{10} 5017600\times60^{-4}$,
$d^{2}/b^{2} = 01~23~13~46~40 =_{10} 17977600\times60^{-4}$.
\\\\\indent We note by taking the \emph{talkiltum} of the short side and the diagonal:
$
a^{2}/b^{2} = \tilde{a}^{2} \times 60^{-n},
$
$
d^{2}/b^{2} = \tilde{d}^{2} \times 60^{-n},
$
therefore the integer\footnote{Or simply the talkiltums if we apply it the place value notation.} 
\emph{talkiltums} (or our auxiliary variables) $\tilde{d}^{2}$, $\tilde{a}^{2}$ satisfy:
\begin{equation}
\tilde{d}^{2} = \tilde{a}^{2} + 60^{n}.
\end{equation}
At this point we realize by considering the number of sexagesits in a given sexagesimal number 
that in order to generate the 4th column entries is easier to use the short side $a$ in calculations 
instead of the diagonal $d$, therefore we will calculate using $a$ but at the end we simply put the 
sexagesimal $\one = 01$ at the beginning of the sexagesimal numbers.
\subsection{Plimpton Machine Accuracy $\epsilon_{m}$.}
As we can see in
\\\citep{Quinonez2020-0}.
If we take the 4th column in P322 as $d^{2}/b^{2}$ then it looks like a normalized mantissa in sexagesimal base 
instead of binary as in
\citep{IEEE754-1985},
\citep{IEEE754-2008},
\citep{IEEE754-2019};
on the other hand if the 4th column is $a^{2}/b^{2}$ then we have also a mantissa.
Let $x\in \mathbb{R}$, in modern computing machines floating point numbers are represented 
\begin{equation*}
x = s \cdot M \cdot B^{e-E},
\end{equation*}
where $B$ is the base, $s$ a bit left for the sign, $e$ the exponent and $E$ the bias in the exponent.
The mantissa $M$ can be taken normalized with $M \in [\frac{1}{B}, 1)$.
In both cases normalized or not, the P322 sexagesimal mantissa does not use exponent and does not use
the $s$ for the sign because they are dealing with positive integers. Therefore we have
\begin{equation}
x = M,
\end{equation}
thus put it in sexagesimal base
\begin{equation}
M = 
\begin{cases}
m_{1} \mathsf{10^{-1}} +  \ldots + m_{8} \mathsf{10^{-8}}, \quad\text{not normalized,}\\
m_{1} \mathsf{10^{0}} + \ldots + m_{9} \mathsf{10^{-8}}, \quad\text{normalized, } m_{1}=1.\\
\end{cases}
\end{equation}
Here we hold the duality of the 4th column representing $a^{2}/b^{2}$ or $d^{2}/b^{2}$,
in both cases the P322 machine accuracy is:
\begin{equation}
\epsilon_{m} = \mathsf{10^{-8}} \approx_{10} 5.95 \times 10^{-15}.
\end{equation}
Columns 2 and 3 (and also 4) have sexagesimal integers in the usual way:
\begin{eqnarray}
x &=& x_{u-1} \mathsf{10^{u-1}} + \ldots + x_{1} \mathsf{10^{1}} + x_{0}, \quad\text{ pure integer},\\
x &=& x_{u-1} \mathsf{10^{-1}} + \ldots + x_{1} \mathsf{10^{-u+1}} + x_{0}\mathsf{10^{-u}}, \quad\text{pure float},
\end{eqnarray}
here the real $x$ is written with $u$ sexagesits as 
$x_{u-1} \ldots x_{1} x_{0}$.
But we can also consider these numbers as mantissas of real numbers, as in the case for the 4th column. 
The converse is alse true. This is one of the advantages of the Old Babylonian Place Value Notation, we could treat
every real number with finite sexagesits as integer!
%
%
%
%
%
\section{The Algorithm of Factor 12.}
The original P322 does not show a column with the values of the side $b$. 
From experience working with data you ommit to show a piece of information in a published (carved) work if either
a) it is not known, or b) it is obvius, perhaps repeated or maybe with the same pattern for all rows, then we reserve the space,
or c) it is deductible from the information already presented. In this opportunity we have all the three cases in one! 
Let us rewrite Table~\ref{t:Plimpton322} with arabic numbers in sexagesimal base
but including the $b$ values and our goal is to construct the table from right to the left.
\begin{table}[h!]
\centering
\begin{tabular}{llllc}
\hline
$\frac{d^{2}}{b^{2}} \big| \frac{d^{2}}{b^{2}} - 1 = \frac{a^{2}}{b^{2}}$ & $a$ & $d$ & $b=12\times Q$ & as itself \\
\hline
59~00~15                & 01~59     & 02~49    & $12\times 10$ & 01 \\
56~56~58~14~56~15       & 56~07     & 01~20~25 & $12\times 04~48$ & 02 \\
55~07~41~15~33~45       & 01~16~41  & 01~50~49 & $12\times 06~40$ & 03 \\
53~10~29~32~52~16       & 03~31~49  & 05~09~01 & $12\times 18~45$ & 04 \\
58~54~01~40             & 01~05     & 01~37    & $12\times 06$ & 05 \\
47~06~41~40             & 05~19     & 08~01    & $12\times 30$ & 06 \\
43~11~56~28~26~40       & 38~11     & 59~01    & $12\times 03~45$ & 07 \\
41~33~45~14~03~45       & 13~19     & 20~49    & $12\times 01~20$ & 08  \\
38~33~36~36             & 08~01     & 12~49    & $12\times 50$ & 09  \\
35~10~02~28~27~24~26~40 & 01~22~41  & 02~16~01 & $12\times 09~00$ & 10 \\
33~45                   & 45        & 01~15    & $12\times 05$ & 11  \\
29~21~54~02~15          & 27~59     & 48~49    & $12\times 03~20$ & 12  \\
27~00~03~45             & 02~41     & 04~49    & $12\times 20$ & 13  \\
25~48~59~35~06~40       & 29~31     & 53~49    & $12\times 03~45$ & 14  \\
23~13~46~40             & 28~00     & 53~00    & $12\times 03~45$ & 15 \\
\hline
\end{tabular}
\caption{Plimpton 322 tablet in sexagesimal basis. We note the $b$ values are scaled by a factor of 12.}
\label{tab:P2}
\end{table}
We see from Table~\ref{tab:P2} that all $b$ values are of the form 
\begin{equation}
b=12\times Q,
\label{e:beq12Q}
\end{equation}
where $Q$ is the scale factor or scale generator of $b$.
%
\subsection{Developing the Algorithm.}
In this section we shall exploit the fact that $12 \divides b$.
For instance we have the diagonal rule $d^{2}=a^{2}+b^{2}$ with $a<b<d$, now we put the emphasis in the large side $b$ and get
\begin{equation}
(d-a)(d+a)=b^{2},
\label{e:b2}
\end{equation}
as a binary multiplicative partition for $b^{2}$. 
If we call
\begin{eqnarray}
y &=& (d+a), \label{e:y} \\
x &=& (d-a), \label{e:x} \\
Y &=& y/Q, \label{e:Y}\\
X &=& x/Q, \label{e:X}
\end{eqnarray}
then inserting (\ref{e:y}) and (\ref{e:x}) into (\ref{e:b2})
\begin{equation}
xy = b^{2};
\label{e:xy}
\end{equation}
and inserting (\ref{e:beq12Q}), (\ref{e:Y}), (\ref{e:X}) into (\ref{e:xy})
\begin{equation}
X Y = 12^{2}.
\label{e:XYeq12^2}
\end{equation}
Now we can solve either the original (\ref{e:xy}) or (\ref{e:XYeq12^2}) and then find appropriate $Q$. 
The system (\ref{e:y}) and (\ref{e:x}) can be solved:
\begin{eqnarray}
d &=& (y+x)/2 = (12^{2}Q^{2}/x + x)/2,  \label{e:d} \\
a &=& (y-x)/2 = (12^{2}Q^{2}/x - x)/2,, \label{e:a} 
\end{eqnarray}
or equivalently, dividing (\ref{e:d}) and (\ref{e:a}) by the scale factor $Q$:
\begin{eqnarray}
D &=& (Y+X)/2 = (12^{2}/X + X)/2,  \label{e:D} \\
A &=& (Y-X)/2 = (12^{2}/X - X)/2. \label{e:A} 
\end{eqnarray}
%
Then the two right triangles $(A, B=12, D)$ and $(a, b, d)$ are similar, therefore
\begin{equation}
\frac{d^{2} }{b^{2}} = \frac{ (Y+X)^{2} }{2^{2} \cdot 12^{2} } = \frac{ (12/X + X/12)^{2} }{2^{2}} = \frac{D^{2}}{12^{2}} =  5 \times D  \times 5 \times D \times 60^{-2}, \label{e:d2/b2} 
\end{equation} 
\begin{equation}
\frac{a^{2} }{b^{2}} = \frac{ (Y-X)^{2} }{2^{2} \cdot 12^{2} } = \frac{ (12/X - X/12)^{2} }{2^{2}} = \frac{A^{2}}{12^{2}} =  5 \times A \times 5 \times A \times 60^{-2}. \label{e:a2/b2} 
\end{equation}
Any of equations (\ref{e:d2/b2}) or (\ref{e:a2/b2}) can be used to calculate the numerical values for the 
most to the left column. After that we need to find appropriate values $Q$ in order to get the values for 
$d$ and $a$ from $D$ and $A$.
This can be done by multiplying each side $A, B, D$ times $60^{n}$ with $n$ a positive integer such that
the resulting three values $(A', B', D')$ be integers, then we take the greatest common divisor among the the three sides
and we get $(a, b, d)=(A', B', D')/\gcd(A', B', D')$. Therefore 
\begin{eqnarray}
Q &=& 60^{n}/\gcd( A \times 60^{n}, 12 \times 60^{n}, D \times 60^{n} ), \label{e:Q} \\
d &=& D Q, \label{e:d=DQ}\\
a &=& A Q. \label{e:a=AQ}
\end{eqnarray}
Finally if we are dealing we $b=12Q$ and triple $(a,b,d)$ with $0 < \theta < \pi/4$ then we have $a<b<d$, but
if $\pi/4 < \theta < \pi/2$ then we have $b<a<d$. 
\subsection{The Case $Q=1$.}
The general case is well depicted until now via problem (\ref{e:XYeq12^2}) until the end of the last subsection just above solving for $D, A$ and $A^{2}/12^{2}$.
We will see this method, as all its predecessors that use scale factors in the three sides at the end of their algorithms, 
are not the best in the sense of the number of operations.
The most tedious part is to find the GCD and then the $Q$ and after that to find $d,a$. 
If we use the same scale factor $Q$ that Old Babylonians took,  then the following factors are easily obtained.
Bruin-Robson algorithm 
used\footnote{For example, if we take her approach from our point of view as using bundle factor $M=1$, or $M=60$, or $M=60^{2}$ etc. 
then she searched for values $Q_{1}$, or $Q_{60}$, or $Q_{3600}$, etc.} 
$M=60^{n}$ by using what is called the standard multipliers.
Here we can cheat\footnote{By using factors that we know give the P322 solutions.} exploting the factorization in (\ref{e:beq12Q}), 
\\\\\indent The following $Q$, $\bar{Q}$, $b$, $\bar{b}$, are factors that can be used in order to get the P322 solutions in the frame of the twelve factor algorithm:
\begin{center}
$Q=_{10}5$, $\bar{Q}=_{10}12$; $b=_{10}60$, $\bar{b} =_{10} 1$.
\\ $Q=_{10}6$, $\bar{Q}=_{10}10$; $b=_{10}72$, $\bar{b}=_{10} 50$.
\\ $Q=_{10}10$, $\bar{Q}=_{10}6$; $b=_{10}120$, $\bar{b} =_{10} 30$.
\\ $Q=_{10}20$, $\bar{Q}=_{10}3$; $b=_{10}240$, $\bar{b} =_{10} 15$.
\\ $Q=_{10}50$, $\bar{Q}=_{10}72$; $b=_{10}600$, $\bar{b} =_{10} 6$.
\\ $Q=_{10}80$, $\bar{Q}=_{10}45$; $b=_{10}960$, $\bar{b} =_{10} 225$.
\\ $Q=_{10}200$, $\bar{Q}=_{10}18$; $b=_{10}2400$, $\bar{b} =_{10} 90$.
\\ $Q=_{10}225$, $\bar{Q}=_{10}16$; $b=_{10}2700$, $\bar{b} =_{10} 80$.
\\ $Q=_{10}288$, $\bar{Q}=_{10}750$; $b=_{10}3456$, $\bar{b} =_{10} 3750$.
\\ $Q=_{10}540$, $\bar{Q}=_{10}400$; $b=_{10}6480$, $\bar{b} =_{10} 2000$.
\\ $Q=_{10}1125$, $\bar{Q}=_{10}192$; $b=_{10}13500$, $\bar{b}=_{10} 960$.
\end{center}
Then we would have the $Q$ values, but instead of that 
in this case $Q=1$ we shall use a nice factor $01.~12$ that is simply $1.2$ or $120\%$ in the following minialgorithm.
\\\\\indent The true is we shall learn at the end of this work, see Figure~\ref{f:histoallQ}, that there are a lot of $Q$ values that give us
a lot of solutions.
Then surges the natural question, why those P322 $Q$? 
As Robson explained those $Q$ are quantities with a reciprocal $\overline{Q}$; we agree
but also there are many others that have a reciprocal that are not presented in P322, then must exist a criterion to choose the fifteen solutions.
\\\\\indent Besides the tedious part mentioned earlier, with the algorithm of factor 12, once found $Y$, $X$ we can find the 4th column values 
in a straightforward way as $A^{2}/12^{2}$.
Studying the case $Q=1$ will give us information about the bounds of the variables $X,Y$ and then we can scale these bounds.
\label{s:Q=1}
\subsubsection{Theoretical Auxiliary Partition or Multiplication Tables for $12^{2}$.}
A simple way to obtain several pairs $X,Y$ such that $XY=02~24$ at once, is by doing a simple modification of the standard multiplication 
Table \ref{tab:std} 
then we multiply the $n$ or $\overline{n}$ columns by 
$12/5=12\times 12\times 01~00^{-1}=02~24\times 01~00^{-1}$.
With respect to this step, there is also the possibility to multiply $n$ (or $\overline{n}$) by $6/5$ and its reciprocal $\overline{n}$ (or $n$) by $2$, 
etc.
Thus the evidence support us because a lot of multiplication tables for 
$6/5=01.~12$ had been found and even applications to interpolation
\citep{Knuth1972}, \citep{AMS-MT}.
\begin{table}[h!]
\centering
\begin{tabular}{llllll}
\hline
$n$ & $\bar{n}$ & $n$ & $\bar{n}$ & $n$ & $\bar{n}$ \\
\hline
02 & 30      & 16 & 03.~45    & 45     & 01.~20 \\
03 & 20      & 18 & 03.~20    & 48     & 01.~15 \\    
04 & 15      & 20 & 03        & 50     & 01.~12 \\    
05 & 12      & 24 & 02.~30    & 54     & 01.~06~40 \\    
06 & 10      & 25 & 02.~24    & 01.~00 & 01~00 \\    
08 & 07.~30  & 27 & 02.~13~20 & 01.~04 & 56.~15 \\    
09 & 06.~40  & 30 & 02        & 01.~12 & 50 \\    
10 & 06      & 32 & 01.~52~30 & 01.~15 & 48 \\    
12 & 05      & 36 & 01.~40    & 01.~20 & 45 \\    
15 & 04      & 40 & 01.~30    & 01.~21 & 44.~26~40 \\    
\hline
\end{tabular}
\caption{Standard Multiplication Table of Reciprocals.
This version is taken from Sh\/{o}yen collection MS3874 in \citep{AMS-MT} this part is due to Friberg. }
\label{tab:std}
\end{table}
Therefore let $U$ and $V$ be a binary multiplicative partition for $02.~24$, 
\begin{equation}
UV = 02.~24
\label{eq:UV}
\end{equation}
we can find the binary multiplicative partition table for $02~24$ from the
standard multiplication table of reciprocals Table~\ref{tab:std} by
requiring the direct product satisfies
\begin{equation}
n U \times \bar{n} V=02~24.
\label{eq:nunv}
\end{equation}
Then the result (\ref{eq:nunv}) is got
multiplying its parts
by $02~24$ times the base $01~00$ powered to $-1$.
This recursive relation 
is provided by the factorizing properties of the 
sexagesimal base $01~00=60_{10}$ and by the right choice of the factor 12.
\\\\\indent There is a lot of ways to start to find right triangle triples such that give the same angles as the triples in P322, but lets choice:
$U=1$, $V=12/5=2\times 01.~12$, then pairs $(n,\bar{n})$ in Table~\ref{tab:std}:
$(05,\: 12)$, $(06,\: 10)$, $(06~40,\: 09)$, $(09,\: 06.~40)$, $(25, \: 02.~24)$, $(27, \: 02.~13~20)$,
$(32,\: 01.~52~30)$, $(40,\: 01.~30)$
become in pairs $(nU, \bar{n}V)$:
$(05, \: 28.~48)$, $(06, \: 24)$, $(06~40, \: 21.~36)$, $(09,\: 16)$, $(25, \: 05.~45~36)$, $(27, \: 05.~20)$, 
$(32,\: 04.~30)$, $(40,\: 03.~36)$.
Here the pairs: $(05, \: 28.~48)$, $(06, \: 24)$, except this $(06~40, \: 21.~36)$, $(25, \: 05.~45~36)$, $(27, \: 05.~20)$ 
all have at least one element at order $60^{0}$ 
and all five without exception will give us
P322 scaled triangles. 
So we can say we have found four duples $X$, $Y$ as solutions with one element at order zero.
Then we can search for constants $k$ and we are looking to find new duples 
\begin{equation}
k \star (nU, \bar{n}V) \equiv (knU, \bar{n}V/k),
\end{equation}
and then we have
$$(01.~00~45) \star (05,\:\: 28.~48) = (05.~03~45,\:\: 28.~26~40),$$
$$(01.~01~26~24) \star (05,\:\: 28.~48) = (05.~07~12,\:\: 28.~07~30),$$
$$(01.~12)^{-2} \star (40,\:\: 01.~30) = (27.~46~40,\:\: 05.~11~02~24),$$
$$(01.~12)^{-1} \star (32,\:\: 04.~30) = (26.~40,\:\: 05.~24),$$
$$(01.~12)^{-1} \star (06.~40,\:\: 21.~36) = (05.~33~20,\:\: 25.~55~12),$$
$$(01.~07~30) \star (05,\:\: 28.~48) = (05.~37~30,\:\: 25.~36),$$
$$(01.~11~06~40) \star (05,\:\: 28.~48) = (05.~55~33~20,\:\: 24.~18),$$
$$(01.~12)^{-2} \star (09,\:\: 16) = (06.~15,\:\: 23.~02~24),$$
$$(01.~12)^{-1} \star (27,\:\: 05.~20) = (22.~30,\:\: 06.~24),$$
$$(01.~12)^{-2} \star (32,\:\: 04.~30) = (22.~13~20,\:\: 06.~28~48).$$

So we have found for the Plimpton 322 generators in the 12-factor algorithm $X \in [05, 06.~40]$, $Y \in [21.~36, 28.~48]$ where they take values according to the 
exact sexagesimal trigonometry for a right triangle with large side with length value 12.
\\\\However there is another property that we have not exploited, if $X,Y$ is a pair solution $XY=02~24$
and we want to find a new solution pair $J,K$ such that $JK=02~24$, with $J=X+E$ for a small known $E$, and $K=Y-S$, we get
\begin{equation}
S = \frac{EY}{J},
\label{eq:S}
\end{equation}
or for a known $S$ and unknown $E$ we get,
\begin{equation}
E = \frac{SX}{K}.
\label{eq:E}
\end{equation}
Finally we gathered the results in Table~\ref{tab:stdsol} and note the pattern in the factors with respect to the solutions
with $X=5$ and $X=6$.
\begin{table}[h!]
\centering
\begin{tabular}{llllllc}
\hline
$X/5$ & $X/6$ & $X$ & $Y$ & $A=\frac{(Y-X)}{2}$  & $D=\frac{(Y+X)}{2}$ & Place\\
\hline
01 & - & 05 & 28.~48 & 11.~54 & 16.~54  & 1 \\
01.~00~45 & - & 05.~03~45 & 28.~26~40  & 11.~41~27~30 & 16.~45~02~30  & 2 \\
01.~01~26~24 & - & 05.~07~12 & 28.~07~30  & 11.~30~09 & 16.~37~21 & 3 \\
01.~02~12~28~48 & - & 05.~11~02~24 & 27.~46~40 & 11.~19~49~18 & 16.~28~54~12 & 4 \\
01.~04 & - & 05.~20 & 27 & 10.~50 & 16.~10 & 5 \\
01.~04~48 & - & 05.~24 & 26.~40 & 10.~38 & 16.~02 & 6 \\
01.~06~40 & - & 05.~33~20 & 25.~55~12 & 10.~10~51 & 15.~49~21 & 7 \\
01.~07~30 & - & 05.~37~30 & 25.~36 & 09.~59~20 & 08.~18~20 & 8 \\
01.~09~07~12 & - & 05.~45~36 & 25 & 09.~37~12 & 15.~22~48 & 9 \\
01.~11~06~40 & - & 05.~55~33~20 & 24.~18 & 09.~11~13~20 & 15.~06~46~40 & 10 \\
01.~12 & 01 & 06 & 24 & 09 & 15  & 11\\
01.~15 & 01.~02~30 & 06.~15 & 23.~02~24 & 08.~53~42 & 14.~38~42 & 12\\
01.~16~48 & 01.~04 & 06.~24 & 22.~30 & 08.~03 & 14.~27 & 13 \\
01.~17~45~36 & 01.~04~48 & 06.~28~48 & 22.~13~20 & 07.~52~16 & 19.~21~04 & 14 \\
01.~20 & 01.~06~40 & 06.~40 & 21.~36 & 07.~28 & 14.~08 & 15  \\
\hline
01.~21 & 01.~07~30 & 06.~45 & 21.~10 & 07.~27~30 & 14.~12~30 & none  \\
\hline
\end{tabular}
\caption{Binary multiplicative partition table for $12^{2}$.}
\label{tab:stdsol}
\end{table}
\subsection{Solutions in other Numerical Systems.} 
\label{s:numericalsys}
We shall require the same conditions in order to get $(a/b)^{2}$ with accuracy of approx. 
$$60^{-8} \approx 30^{-10} \approx 16^{-12} \approx 10^{-15} \approx 8^{-16} \approx 2^{-48},$$ 
this means to require 8, 12, 15, 16, and 48 places of significance for the bases sexagesimal, hexadecimal, decimal, octal and binary respectively;
also that the floating point representation for the pair numbers $X_{M}$, $Y_{M}$ in the given base ends in that significance position or before,
i.e. the pair must have a finite floating point representation with the bounding accuracy.
Then we have to be able to find an appropriate mak\c{s}arum $M$ for the chosen numerical base.
To find these solutions is a very time consuming task and we can get the solutions with help of a computer, but if we restrict ourselves
to use only the methods that Old Babylonians get used to do, the hard task is to find new binary multiplicative partition tables with all the required 
places after the point radix, but (\ref{eq:S}) and (\ref{eq:E}) would help in order to get new solutions. 
Certainly there are exists another sets of solutions with not necessarily the same proportioned right triangles or angles.
We shall see some of them in Section \ref{s:solutions}.
\subsection{Out of the P322 Range or Another Angle Solutions.}
\label{s:p322}
Lets $z_{2}$ and $z_{1}$ be small numbers in $\mathbb{R}$ with finite representation on base $B$ up to order $B^{-4}$,
we can write:
\begin{eqnarray*}
z_{1} &=& a_{1} + b_{1} \times B^{-1} + c_{1} \times B^{-2} + d_{1} \times B^{-3} + e_{1} \times B^{-4}   \\
z_{2} &=& a_{2} + b_{2} \times B^{-1} + c_{2} \times B^{-2} + d_{2} \times B^{-3} + e_{2} \times B^{-4}   
\end{eqnarray*}
then 
\begin{eqnarray}
z_{2} \pm z_{1} &=& (a_{2} \pm a_{1} ) + (b_{2} \pm b_{1}) \times B^{-1} + (c_{2} \pm c_{1}) \times B^{-2}\\
& & + (d_{2} \pm d_{1}) \times B^{-3} + (e_{2}\pm e_{1}) \times B^{-4}. \nonumber
\label{eq:}
\end{eqnarray}
Then the product of two different sexagesimal numbers with order up to $B^{-4}$ is:
\begin{eqnarray}
z_{1} z_{2} &=& a_{1}a_{2} + (a_{1}b_{2}+b_{1}a_{2}) \times B^{-1} \\
 & & + (b_{1}b_{2} + a_{1}c_{2}+c_{1}a_{2}) \times B^{-2} \nonumber\\
 & &  + (a_{1}d_{2}+d_{1}a_{2}+b_{1}c_{2}+c_{1}b_{2}) \times B^{-3}   \nonumber\\
 & &  + (c_{1}c_{2}+b_{1}d_{2}+d_{1}b_{2}+a_{1}e_{2}+e_{1}a_{2}) \times B^{-4}   \nonumber\\
 & &  + (c_{1}d_{2}+d_{1}c_{2}+b_{1}e_{2}+e_{1}b_{2}) \times B^{-5}   \nonumber\\
 & &  + (d_{1}d_{2}+c_{1}e_{2}+e_{1}c_{2}) \times B^{-6}  \nonumber\\
 & &  + (d_{1}e_{2}+e_{1}d_{2}) \times B^{-7}  + (e_{1}e_{2}) \times B^{-8},   \nonumber
\label{eq:}
\end{eqnarray}
thus we get
\begin{eqnarray}
(z_{2} \pm z_{1})^{2} &=& (a_{2} \pm a_{1})^{2} + [ 2 (a_{2} \pm a_{1}) (b_{2} \pm b_{1}) ] \times B^{-1} \label{e:zz} \\
 & & + [ (b_{2} \pm b_{1})^{2} + 2 (a_{2}\pm a_{1})(c_{2} \pm c_{1}) ] \times B^{-2} \nonumber\\
 & &  + [2(a_{2} \pm a_{1})(d_{2} \pm d_{1})  + 2(b_{2} \pm b_{1})(c_{2} \pm c_{1}) ] \times B^{-3}   \nonumber\\
 & &  + [(c_{2} \pm c_{1})^{2} + 2(b_{2} \pm b_{1})(d_{2} \pm d_{1}) + 2(a_{2} \pm a_{1})(e_{2} \pm e_{1} ) ] \times B^{-4}   \nonumber\\
 & &  + [2(c_{2} \pm c_{1})(d_{2} \pm d_{1}) +2 (b_{2} \pm b_{1})(e_{2} \pm e_{1}) ] \times B^{-5}   \nonumber\\
 & &  + [(d_{2} \pm d_{1})^{2} + 2(c_{2} \pm c_{1})(e_{2} \pm e_{1}) ] \times B^{-6}  \nonumber\\
 & & + [2(d_{2} \pm d_{1})(e_{2} \pm e_{1} )] \times B^{-7}  + (e_{2} \pm e_{1})^{2} \times B^{-8},   \nonumber
\label{eq:}
\end{eqnarray}
therefore the product of factor 
$$
\frac{1}{2^{2}}\frac{1}{12^{2}} = 30 \times 30 \times 05 \times 05 \times 01~00^{-4} = 06~15\times 01~00^{-3} 
$$ 
with $(\ref{e:zz})$ and the condition to get P322 solutions with generators $X,Y$ with accuracy up to $01~00^{-3}$,
lead us to conclude: 
\\\\\indent \emph{$z_{2}=Y \in [28.~48, \:\: 21.~36] $ and $z_{1}=X \in [05, \:\: 06.~40]$, both with order up to $01~00^{-3}$, 
is a neccesary but not sufficient condition to get P322 solutions in 4th column with accuracy 
up to $01~00^{-8}$.} 
\\\\Tables~\ref{t:PextA}, \ref{t:PextB} show solutions in this scheme by using (\ref{e:a2/b2}), (\ref{e:Q}), (\ref{e:d=DQ}), (\ref{e:a=AQ})
for interval $\pi/4 > \theta > 0$.
\subsubsection{Application: The Angle of the Great Pyramid of Giza.} 
When applying the 12-factor algorithm for P322
at fourth order in $X, Y$, it is found an interesting solution for the Great Pyramid of Giza 
$X=05.~49~55~12$, $Y=24.~41~28~53~20$, $A=_{10} 9.4296790123456$, $D=_{10} 15.261679012345$, 
$Q = 05~37~30=_{10} 20250$, $b = 01~07~30~00 = 243000$, $a=53~02~31=_{10} 190951$, $d=01~25~50~49=_{10} 309049$,
therefore 
$\theta =_{10} 0.666026302753977$, $\phi =_{10} 0.904770024040919$. The former is the Giza apothem-base angle which in degrees is:
$51.83950380749469^{\circ}$. 
This value lies in the middle of the values coming from the $\pi$ theory $51.854^{\circ}$ and $\phi$ theory $51.827$, see \citep{Bartlett2014}.
\begin{table}
\centering
\begin{tabular}{cllll}
\hline
Place & $X$ & $Y$ & $A=(Y-X)/2$ & $D=(Y+X)/2$   \\
\hline
1*	&	5	&	 28.~48 	&	 11.~54 	&	 16.~54 	\\
2*	&	 05.~03~45 	&	 28.~26~40 	&	 11.~41~27~30 	&	 16.~45~12~30 	\\
3*	&	 05.~07~12 	&	 28.~07~30 	&	11.~30~09	&	16.~37~21	\\
4*	&	 05.~11~02~24 	&	 27.~46~40 	&	11.~17~48~48	&	16.~28~51~12	\\
5	&	 05.~12~30 	&	 27.~38~52~48 	&	11.~13~11~24	&	16.~25~41~24	\\
6*	&	 05.~20 	&	27	&	10.~50	&	16.~10	\\
7*	&	 05.~24 	&	 26.~40 	&	10.~38	&	16.~02	\\
8*	&	 05.~33~20 	&	 25.~55~12 	&	10.~10~56	&	15.~44~16	\\
9*	&	 05.~37~30 	&	 25.~36 	&	09.~59~15	&	15.~36~45	\\
10	&	 05.~41~20 	&	 25.~18~45 	&	09.~48~42~30	&	15.~30~02~30	\\
11*	&	 05.~45~36 	&	25	&	09.~37~12	&	15.~22~48	\\
12*	&	 05.~55~33~20 	&	 24.~18 	&	09.~11~13~20	&	15.~06~46~40	\\
13*	&	6	&	24	&	9	&	15	\\
14	&	 06.~04~30 	&	 23.~42~13~20 	&	08.~48~51~40	&	14.~53~21~40	\\
15	&	 06.~08~38~24 	&	 23.~26~15 	&	08.~38~48~18	&	14.~47~26~42	\\
16*	&	 06.~15 	&	 23.~02~24 	&	08.~23~42	&	14.~38~42	\\
17*	&	 06.~24 	&	 22.~30 	&	08.~03	&	14.~27	\\
18*	&	 06.~28~48 	&	 22.~13~20 	&	07.~52~16	&	14.~21~04	\\
19*	&	 06.~40 	&	 21.~36 	&	07.~28	&	14.~08	\\
\hline
20	&	 06.~45 	&	 21.~20 	&	07.~17~30	&	14.~02~30	\\
21	&	 06.~49~36 	&	 21.~05~37~30 	&	07.~08~00~45	&	13.~57~36~45	\\
22	&	 06.~56~40 	&	 20.~44~09~36 	&	06.~53~44~48	&	13.~50~24~48	\\
23	&	 07.~06~40 	&	 20.~15 	&	06.~34~10	&	13.~40~50	\\
24	&	 07.~12 	&	20	&	06.~24	&	13.~36	\\
25	&	 07.~30 	&	 19.~12 	&	05.~51	&	13.~21	\\
26	&	 07.~35~37~30 	&	 18.~57~46~40 	&	05.~41~04~35	&	13.~16~42~05	\\
27	&	 07.~46~33~36 	&	 18.~31~06~40 	&	05.~22~16~32	&	13.~08~50~08	\\
28	&	8	&	18	&	5	&	13	\\
29	&	 08.~06 	&	 17.~46~40 	&	04.~50~20	&	12.~56~20	\\
30	&	 08.~26~15 	&	 17.~04 	&	04.~18~52~30	&	12.~45~07~30	\\
31	&	 08.~32 	&	 16.~52~30 	&	04.~10~15	&	12.~42~15	\\
32	&	 08.~53~20 	&	 16.~12 	&	03.~39~20	&	12.~32~40	\\
33	&	9	&	16	&	03.~30	&	12.~30	\\
34	&	 09.~15~33~20 	&	 15.~33~07~12 	&	03.~08~46~56	&	12.~24~20~16	\\
35	&	 09.~22~30 	&	 15.~21~36 	&	02.~59~33	&	12.~22~03	\\
36	&	 09.~28~53~20 	&	 15.~11~15 	&	02.~51~10~50	&	12.~20~04~10	\\
37	&	 09.~36 	&	15	&	02.~42	&	12.~18	\\
38	&	 09.~43~12 	&	 14.~48~53~20 	&	02.~32~50~40	&	12.~16~02~40	\\
39	&	10	&	 14.~24 	&	02.~12	&	12.~12	\\
40	&	 10.~07~30 	&	 14.~13~20 	&	02.~02~55	&	12.~10~25	\\
41	&	 10.~14~24 	&	 14.~03~45 	&	01.~54~40~30	&	12.~09~04~30	\\
42	&	 10.~22~04~48 	&	 13.~53~20 	&	01.~45~37~36	&	12.~07~42~24	\\
43	&	 10.~25 	&	 13.~49~26~24 	&	01.~42~13~12	&	12.~07~13~12	\\
44	&	 10.~40 	&	 13.~30 	&	01.~25	&	12.~05	\\
45	&	 10.~48 	&	 13.~20 	&	01.~16	&	12.~04	\\
46	&	 11.~06~40 	&	 12.~57~36 	&	00.~55~28	&	12.~02~08	\\
47	&	 11.~15 	&	 12.~48 	&	00.~46~30	&	12.~01~30	\\
48	&	 11.~22~40 	&	 12.~39~22~30 	&	00.~38~21~15	&	12.~01~1~15	\\
49	&	 11.~31~12 	&	 12.~30 	&	00.~29~24	&	12.~00~36	\\
50	&	 11.~43~07~30 	&	 12.~17~16~48 	&	00.~17~04~39	&	12.~00~12~09	\\
51	&	 11.~51~06~40 	&	 12.~09 	&	00.~08~56~40	&	12.~00~03~20	\\
\hline
\end{tabular}
\caption{Solutions for interval $\pi/4 > \theta > 0$ up to order $60^{-3}$ in $X, Y$.}
\label{t:PextA}
\end{table} 

\begin{table}
\centering
\begin{tabular}{cllll}
\hline
Place & $b$ & $d$ & $a$ & $a^{2}/b^{2} = d^{2}/b^{2} - 1 $   \\
\hline
1*	&	120	&	169	&	119	&	0.983402777777778	\\
2*	&	3456	&	4825	&	3367	&	0.949158552088692	\\
3*	&	4800	&	6649	&	4601	&	0.918802126736111	\\
4*	&	13500	&	18541	&	12709	&	0.886247906721536	\\
5	&	72000	&	98569	&	67319	&	0.874199027970679	\\
6*	&	72	&	97	&	65	&	0.815007716049383	\\
7*	&	360	&	481	&	319	&	0.785192901234568	\\
8*	&	2700	&	3541	&	2291	&	0.719983676268862	\\
9*	&	960	&	1249	&	799	&	0.692709418402778	\\
10	&	17280	&	22321	&	14129	&	0.668552279583548	\\
11*	&	600	&	769	&	481	&	0.642669444444444	\\
12*	&	6480	&	8161	&	4961	&	0.586122566110349	\\
13*	&	4	&	5	&	3	&	0.5625	\\
14	&	25920	&	32161	&	19039	&	0.539533323390584	\\
15	&	16000	&	19721	&	11529	&	0.51921031640625	\\
16*	&	2400	&	2929	&	1679	&	0.489416840277778	\\
17*	&	240	&	289	&	161	&	0.450017361111111	\\
18*	&	2700	&	3229	&	1771	&	0.430238820301783	\\
19*	&	90	&	106	&	56	&	0.38716049382716	\\
\hline
20	&	288	&	337	&	175	&	0.369225019290123	\\
21	&	57600	&	67009	&	34241	&	0.353384838264371	\\
22	&	54000	&	62281	&	31031	&	0.330220494170096	\\
23	&	864	&	985	&	473	&	0.299705557698903	\\
24	&	30	&	34	&	16	&	0.284444444444444	\\
25	&	80	&	89	&	39	&	0.23765625	\\
26	&	20736	&	22945	&	9823	&	0.224407992009037	\\
27	&	20250	&	22186	&	9064	&	0.200350180765127	\\
28	&	12	&	13	&	5	&	0.173611111111111	\\
29	&	2160	&	2329	&	871	&	0.162603094993141	\\
30	&	5760	&	6121	&	2071	&	0.129275203751929	\\
31	&	2880	&	3049	&	1001	&	0.12080451871142	\\
32	&	1080	&	1129	&	329	&	0.092799211248285	\\
33	&	24	&	25	&	7	&	0.085069444444445	\\
34	&	40500	&	41869	&	10619	&	0.068747545191282	\\
35	&	1600	&	1649	&	399	&	0.062187890625	\\
36	&	10368	&	10657	&	2465	&	0.056525429398744	\\
37	&	40	&	41	&	9	&	0.050625	\\
38	&	16200	&	16561	&	3439	&	0.045064475689682	\\
39	&	60	&	61	&	11	&	0.033611111111111	\\
40	&	1728	&	1753	&	295	&	0.029144496420611	\\
41	&	9600	&	9721	&	1529	&	0.025367198350695	\\
42	&	27000	&	27289	&	3961	&	0.021521976680384	\\
43	&	115200	&	116161	&	14911	&	0.0167536169735	\\
44	&	144	&	145	&	17	&	0.013937114197531	\\
45	&	180	&	181	&	19	&	0.011141975308642	\\
46	&	1350	&	1354	&	104	&	0.005934705075446	\\
47	&	480	&	481	&	31	&	0.004171006944444	\\
48	&	34560	&	34609	&	1841	&	0.002837658373253	\\
49	&	1200	&	1201	&	49	&	0.001667361111111	\\
50	&	32000	&	32009	&	759	&	0.000562579101563	\\
51	&	12960	&	12961	&	161	&	0.000154326941396	\\
\hline
\end{tabular}
\caption{Solutions in base 10 for interval $1 > a^{2}/b^{2} > 0$ up to order $60^{-3}$ in $X, Y$.}
\label{t:PextB}
\end{table} 
\subsection{Algorithms of the Form $b=M \, Q_{M}$.}
Taking the 15th row {\it a la} Robson, then after we get $b=3 \, Q_{3}$.
If we do not consider the fifteen entry then the result is $b=12 \, Q$;
but if we additionally do not consider the second and fifth entries the result is $b=60 \, Q_{60}$.
Therefore we can get another maksarums for another set of P322 solutions, in brief we discuss some of them.
Returning to the 15th row if we consider this entry by using the Old Babylonian place value notation as in Table~\ref{tab:P2}
we recover again $b=12\,Q$.
\subsubsection{$M=3$.} 
In this theory $b=3 \, Q_{3} = 12 \, Q$ and $X_{3}Y_{3}=3^{2}$, therefore $Q_{3}=4 \, Q$ and $X_{3}=X/4$, $Y_{3}=Y/4$.
The range of the monotonically changing generators are:
$X_{3} \in (01.~15, \:\: 07.~12)$,
$Y_{3} \in (01.~40, \:\: 05.~24)$.
\subsubsection{$M=60$.}
In this theory $b=60 \, Q_{60} = 12 \, Q$ and $X_{60}Y_{60}=60^{2}$, therefore $Q_{60} = Q/5$ and $X_{60}=5X$, $Y_{60}=5Y$.
The range of the monotonically changing generators are:
$X_{60} \in (25, \:\: 33~20)$,
$Y_{60} \in (02~24, \:\: 01~48)$.
Or perhaps another wished choice.
%
\begin{table}
\centering
\begin{tabular}{rlll}
\hline
Row & $Q$ & $Q_{3}$ & $Q_{60}$ \\
\hline
1 & 10       & 40    & 02     \\
2 & 04~48    & 19~12 & 57.~36 \\
3 & 06~40    & 26~40 & 01~20  \\
4 & 18~45    & 01~15~00 & 03~45  \\
5 & 06    & 24 & 01.~12  \\
6 & 30    & 02~00 & 06  \\
7 & 03~45    & 15~00 & 45  \\
8 & 01~20    & 05~20 & 16  \\
9 & 50    & 03~20 & 10  \\
10 & 09~00    & 36~00 & 01~48  \\
11 & 05    & 20 & 01  \\
12 & 03~20    & 13~20 & 40  \\
13 & 20    & 03~20 & 04  \\
14 & 03~45    & 15~00 & 45  \\
15 & 03~45    & 15 & 45  \\
\hline
none & 04~48  & 19~12 & 57.~36  \\
\hline
\end{tabular}
\caption{Scale generator values in theories $M=12,3,60$ for P322 entries.}
\label{t:QQQ}
\end{table}
\\\\\indent At this moment we can also note from Tables~\ref{tab:P2},~\ref{t:QQQ} that the $Q$ values are apparently random and not sorted. 
However we shall see in Section~\ref{s:solutions} that in the end we can consider this $Q$ sorted and choiced according to a nice criterion.
\section{Integer Solutions}
\label{s:solutions}
The methodologies applied in Sections~\ref{s:Q=1} and \ref{s:p322} are far beyond to be optimal and
for the same reason also are the methodologies applied in all works up to date.
The most time consuming part is to find the $Q$ values.
In Section~\ref{s:Q=1} we can note the number of sexagesits for $x=X$ and $y=Y$ are greater
than presented in the following tables of solutions where $Q$ has been taken with the same values as in P322
but here we search for all integer solutions from the beginning, 
i.e. we take $x$, $y$, and $Q$ as integers and we aim to find integer values for $(a, b, d)$.
This is the integer solutions scheme and is based on equations (\ref{e:beq12Q}), (\ref{e:xy}), (\ref{e:d}), (\ref{e:a}), and (\ref{e:a2/b2}).
\\When $\pi/4 < \theta < \pi/2$ the side $b$ becomes the small side, and the side $a$ becomes in the new large side.
From the experience gained in Section~\ref{s:Q=1} we know P322 solutions for $X \in [5, 6\frac{40}{60}]$ and
$Y \in [21\frac{36}{60}, 28\frac{48}{60}]$. We could use this fact in order to find sexagesimal solutions by the 
scribe's methods.
Therefore We are going to present all the solutions obtained varying $x \in [2, 12Q-2]$.
The following tables presents the simbol S-n with n a number, this symbol means
$\times 60^{-n}$. There are also row separation lines indicating when the triple has an angle $\theta \in (\pi/6, \pi/4)$.
We write down the integer solutions in decimal base just to show that we can inherit the ancient sexagesimal techniques inside the modern
decimal base in order to get nice integer results with sexagesimal scientific notation.
You can note if we work the solutions from beginning in an ordered way, in terms of the generators $x$, $y$, requiring $xy=b^{2}$ such that 
they be integers, then we  
do not strongly need the $Q$ factors in the arithmetical operations 
to get $y$\footnote{We can think the problem as: find an integer that multiplied by $x$ gives $b^{2}$.}, 
$d=(y+x)/2$, $a=(y-x)/2$, and we need them for $a^{2}/b^{2} = A^{2}/12^{2}$ because we saw in the data the scribe write $a^{2}$ when 
he should write $a$; he also made mistakes by writing down $2a$, this can be explained by noting $2a=y-x$. 
%
\begin{center}
\begin{tabular}{ccccccc}
\hline
$x$ & $y$ & $b=12\cdot 5$ & $a$ & $d$ & $a^{2}$ & $a^{2}/b^{2}$  \\
\hline
2 & 1800 & 60 & 899 & 901 & 808201 & 808201S-2\\
4 & 900 & 60 & 448 & 452 & 200704 & 200704S-2\\
6 & 600 & 60 & 297 & 303 & 88209 & 88209S-2\\
8 & 450 & 60 & 221 & 229 & 48841 & 48841S-2\\
10 & 360 & 60 & 175 & 185 & 30625 & 30625S-2\\
12 & 300 & 60 & 144 & 156 & 20736 & 20736S-2\\
18 & 200 & 60 & 91 & 109 & 8281 & 8281S-2\\
20 & 180 & 60 & 80 & 100 & 6400 & 6400S-2\\
24 & 150 & 60 & 63 & 87 & 3969 & 3969S-2\\
\hline
30 & 120 & 60 & 45 & 75 & 2025 & 2025S-2\\
\hline
36 & 100 & 60 & 32 & 68 & 1024 & 1024S-2\\
40 & 90 & 60 & 25 & 65 & 625 & 625S-2\\
50 & 72 & 60 & 11 & 61 & 121 & 121S-2\\
\hline
\end{tabular}
\end{center}
\begin{center}
\begin{tabular}{ccccccc}
\hline
$x$ & $y$ & $b=12\cdot 6$ & $a$ & $d$ & $a^{2}$ & $a^{2}/b^{2}$  \\
\hline
2 & 2592 & 72 & 1295 & 1297 & 1677025 & 4192562500S-4\\
4 & 1296 & 72 & 646 & 650 & 417316 & 1043290000S-4\\
6 & 864 & 72 & 429 & 435 & 184041 & 460102500S-4\\
8 & 648 & 72 & 320 & 328 & 102400 & 256000000S-4\\
12 & 432 & 72 & 210 & 222 & 44100 & 110250000S-4\\
16 & 324 & 72 & 154 & 170 & 23716 & 59290000S-4\\
18 & 288 & 72 & 135 & 153 & 18225 & 45562500S-4\\
24 & 216 & 72 & 96 & 120 & 9216 & 23040000S-4\\
\hline
32 & 162 & 72 & 65 & 97 & 4225 & 10562500S-4\\
36 & 144 & 72 & 54 & 90 & 2916 & 7290000S-4\\
\hline
48 & 108 & 72 & 30 & 78 & 900 & 2250000S-4\\
54 & 96 & 72 & 21 & 75 & 441 & 1102500S-4\\
\hline
\end{tabular}
\end{center}
\begin{center}
\begin{tabular}{ccccccc}
\hline
$x$ & $y$ & $b=12\cdot 10$ & $a$ & $d$ & $a^{2}$ & $a^{2}/b^{2}$  \\
\hline
2 & 7200 & 120 & 3599 & 3601 & 12952801 & 11657520900S-4\\
4 & 3600 & 120 & 1798 & 1802 & 3232804 & 2909523600S-4\\
6 & 2400 & 120 & 1197 & 1203 & 1432809 & 1289528100S-4\\
8 & 1800 & 120 & 896 & 904 & 802816 & 722534400S-4\\
10 & 1440 & 120 & 715 & 725 & 511225 & 460102500S-4\\
12 & 1200 & 120 & 594 & 606 & 352836 & 317552400S-4\\
16 & 900 & 120 & 442 & 458 & 195364 & 175827600S-4\\
18 & 800 & 120 & 391 & 409 & 152881 & 137592900S-4\\
20 & 720 & 120 & 350 & 370 & 122500 & 110250000S-4\\
24 & 600 & 120 & 288 & 312 & 82944 & 74649600S-4\\
30 & 480 & 120 & 225 & 255 & 50625 & 45562500S-4\\
32 & 450 & 120 & 209 & 241 & 43681 & 39312900S-4\\
36 & 400 & 120 & 182 & 218 & 33124 & 29811600S-4\\
40 & 360 & 120 & 160 & 200 & 25600 & 23040000S-4\\
48 & 300 & 120 & 126 & 174 & 15876 & 14288400S-4\\
\hline
50 & 288 & 120 & 119 & 169 & 14161 & 12744900S-4\\
60 & 240 & 120 & 90 & 150 & 8100 & 7290000S-4\\
\hline
72 & 200 & 120 & 64 & 136 & 4096 & 3686400S-4\\
80 & 180 & 120 & 50 & 130 & 2500 & 2250000S-4\\
90 & 160 & 120 & 35 & 125 & 1225 & 1102500S-4\\
96 & 150 & 120 & 27 & 123 & 729 & 656100S-4\\
100 & 144 & 120 & 22 & 122 & 484 & 435600S-4\\
\hline
\end{tabular}
\end{center}

\begin{center}
$Q=20$
\begin{tabular}{ccccccc}
\hline
$x$ & $y$ & $b=12\cdot 20$ & $a$  & $d$ & $a^{2}$ & $a^{2}/b^{2}$  \\
\hline
2 & 28800 & 240 & 14399 & 14401 & 207331201 & 46649520225S-4\\
4 & 14400 & 240 & 7198 & 7202 & 51811204 & 11657520900S-4\\
6 & 9600 & 240 & 4797 & 4803 & 23011209 & 5177522025S-4\\
8 & 7200 & 240 & 3596 & 3604 & 12931216 & 2909523600S-4\\
10 & 5760 & 240 & 2875 & 2885 & 8265625 & 1859765625S-4\\
12 & 4800 & 240 & 2394 & 2406 & 5731236 & 1289528100S-4\\
16 & 3600 & 240 & 1792 & 1808 & 3211264 & 722534400S-4\\
18 & 3200 & 240 & 1591 & 1609 & 2531281 & 569538225S-4\\
20 & 2880 & 240 & 1430 & 1450 & 2044900 & 460102500S-4\\
24 & 2400 & 240 & 1188 & 1212 & 1411344 & 317552400S-4\\
30 & 1920 & 240 & 945 & 975 & 893025 & 200930625S-4\\
32 & 1800 & 240 & 884 & 916 & 781456 & 175827600S-4\\
36 & 1600 & 240 & 782 & 818 & 611524 & 137592900S-4\\
40 & 1440 & 240 & 700 & 740 & 490000 & 110250000S-4\\
48 & 1200 & 240 & 576 & 624 & 331776 & 74649600S-4\\
50 & 1152 & 240 & 551 & 601 & 303601 & 68310225S-4\\
60 & 960 & 240 & 450 & 510 & 202500 & 45562500S-4\\
64 & 900 & 240 & 418 & 482 & 174724 & 39312900S-4\\
72 & 800 & 240 & 364 & 436 & 132496 & 29811600S-4\\
80 & 720 & 240 & 320 & 400 & 102400 & 23040000S-4\\
90 & 640 & 240 & 275 & 365 & 75625 & 17015625S-4\\
96 & 600 & 240 & 252 & 348 & 63504 & 14288400S-4\\
\hline
100 & 576 & 240 & 238 & 338 & 56644 & 12744900S-4\\
120 & 480 & 240 & 180 & 300 & 32400 & 7290000S-4\\
128 & 450 & 240 & 161 & 289 & 25921 & 5832225S-4\\
\hline
144 & 400 & 240 & 128 & 272 & 16384 & 3686400S-4\\
150 & 384 & 240 & 117 & 267 & 13689 & 3080025S-4\\
160 & 360 & 240 & 100 & 260 & 10000 & 2250000S-4\\
180 & 320 & 240 & 70 & 250 & 4900 & 1102500S-4\\
192 & 300 & 240 & 54 & 246 & 2916 & 656100S-4\\
200 & 288 & 240 & 44 & 244 & 1936 & 435600S-4\\
\hline
\end{tabular}
\end{center}

\begin{center}
\begin{tabular}{ccccccc}
\hline
$x$ & $y$ & $b=12\cdot 30$ & $a$ & $d$ & $a^{2}$ & $a^{2}/b^{2}$  \\
\hline
2 & 64800 & 360 & 32399 & 32401 & 1049695201 & 104969520100S-4\\
4 & 32400 & 360 & 16198 & 16202 & 262375204 & 26237520400S-4\\
6 & 21600 & 360 & 10797 & 10803 & 116575209 & 11657520900S-4\\
8 & 16200 & 360 & 8096 & 8104 & 65545216 & 6554521600S-4\\
10 & 12960 & 360 & 6475 & 6485 & 41925625 & 4192562500S-4\\
12 & 10800 & 360 & 5394 & 5406 & 29095236 & 2909523600S-4\\
16 & 8100 & 360 & 4042 & 4058 & 16337764 & 1633776400S-4\\
18 & 7200 & 360 & 3591 & 3609 & 12895281 & 1289528100S-4\\
20 & 6480 & 360 & 3230 & 3250 & 10432900 & 1043290000S-4\\
24 & 5400 & 360 & 2688 & 2712 & 7225344 & 722534400S-4\\
30 & 4320 & 360 & 2145 & 2175 & 4601025 & 460102500S-4\\
32 & 4050 & 360 & 2009 & 2041 & 4036081 & 403608100S-4\\
36 & 3600 & 360 & 1782 & 1818 & 3175524 & 317552400S-4\\
40 & 3240 & 360 & 1600 & 1640 & 2560000 & 256000000S-4\\
48 & 2700 & 360 & 1326 & 1374 & 1758276 & 175827600S-4\\
50 & 2592 & 360 & 1271 & 1321 & 1615441 & 161544100S-4\\
54 & 2400 & 360 & 1173 & 1227 & 1375929 & 137592900S-4\\
60 & 2160 & 360 & 1050 & 1110 & 1102500 & 110250000S-4\\
72 & 1800 & 360 & 864 & 936 & 746496 & 74649600S-4\\
80 & 1620 & 360 & 770 & 850 & 592900 & 59290000S-4\\
90 & 1440 & 360 & 675 & 765 & 455625 & 45562500S-4\\
96 & 1350 & 360 & 627 & 723 & 393129 & 39312900S-4\\
100 & 1296 & 360 & 598 & 698 & 357604 & 35760400S-4\\
108 & 1200 & 360 & 546 & 654 & 298116 & 29811600S-4\\
120 & 1080 & 360 & 480 & 600 & 230400 & 23040000S-4\\
144 & 900 & 360 & 378 & 522 & 142884 & 14288400S-4\\
\hline
150 & 864 & 360 & 357 & 507 & 127449 & 12744900S-4\\
160 & 810 & 360 & 325 & 485 & 105625 & 10562500S-4\\
162 & 800 & 360 & 319 & 481 & 101761 & 10176100S-4\\
180 & 720 & 360 & 270 & 450 & 72900 & 7290000S-4\\
200 & 648 & 360 & 224 & 424 & 50176 & 5017600S-4\\
\hline
216 & 600 & 360 & 192 & 408 & 36864 & 3686400S-4\\
240 & 540 & 360 & 150 & 390 & 22500 & 2250000S-4\\
270 & 480 & 360 & 105 & 375 & 11025 & 1102500S-4\\
288 & 450 & 360 & 81 & 369 & 6561 & 656100S-4\\
300 & 432 & 360 & 66 & 366 & 4356 & 435600S-4\\
324 & 400 & 360 & 38 & 362 & 1444 & 144400S-4\\
\hline
\end{tabular}
\end{center}

\begin{center}
\begin{tabular}{ccccccc}
\hline
$x$ & $y$ & $b=12\cdot 50$ & $a$  & $d$ & $a^{2}$ & $a^{2}/b^{2}$  \\
\hline
2 & 180000 & 600 & 89999 & 90001 & 8099820001 & 291593520036S-4\\
4 & 90000 & 600 & 44998 & 45002 & 2024820004 & 72893520144S-4\\
6 & 60000 & 600 & 29997 & 30003 & 899820009 & 32393520324S-4\\
8 & 45000 & 600 & 22496 & 22504 & 506070016 & 18218520576S-4\\
10 & 36000 & 600 & 17995 & 18005 & 323820025 & 11657520900S-4\\
12 & 30000 & 600 & 14994 & 15006 & 224820036 & 8093521296S-4\\
16 & 22500 & 600 & 11242 & 11258 & 126382564 & 4549772304S-4\\
18 & 20000 & 600 & 9991 & 10009 & 99820081 & 3593522916S-4\\
20 & 18000 & 600 & 8990 & 9010 & 80820100 & 2909523600S-4\\
24 & 15000 & 600 & 7488 & 7512 & 56070144 & 2018525184S-4\\
30 & 12000 & 600 & 5985 & 6015 & 35820225 & 1289528100S-4\\
32 & 11250 & 600 & 5609 & 5641 & 31460881 & 1132591716S-4\\
36 & 10000 & 600 & 4982 & 5018 & 24820324 & 893531664S-4\\
40 & 9000 & 600 & 4480 & 4520 & 20070400 & 722534400S-4\\
48 & 7500 & 600 & 3726 & 3774 & 13883076 & 499790736S-4\\
50 & 7200 & 600 & 3575 & 3625 & 12780625 & 460102500S-4\\
60 & 6000 & 600 & 2970 & 3030 & 8820900 & 317552400S-4\\
72 & 5000 & 600 & 2464 & 2536 & 6071296 & 218566656S-4\\
80 & 4500 & 600 & 2210 & 2290 & 4884100 & 175827600S-4\\
90 & 4000 & 600 & 1955 & 2045 & 3822025 & 137592900S-4\\
96 & 3750 & 600 & 1827 & 1923 & 3337929 & 120165444S-4\\
100 & 3600 & 600 & 1750 & 1850 & 3062500 & 110250000S-4\\
120 & 3000 & 600 & 1440 & 1560 & 2073600 & 74649600S-4\\
144 & 2500 & 600 & 1178 & 1322 & 1387684 & 49956624S-4\\
150 & 2400 & 600 & 1125 & 1275 & 1265625 & 45562500S-4\\
160 & 2250 & 600 & 1045 & 1205 & 1092025 & 39312900S-4\\
180 & 2000 & 600 & 910 & 1090 & 828100 & 29811600S-4\\
200 & 1800 & 600 & 800 & 1000 & 640000 & 23040000S-4\\
240 & 1500 & 600 & 630 & 870 & 396900 & 14288400S-4\\
\hline
250 & 1440 & 600 & 595 & 845 & 354025 & 12744900S-4\\
288 & 1250 & 600 & 481 & 769 & 231361 & 8328996S-4\\
300 & 1200 & 600 & 450 & 750 & 202500 & 7290000S-4\\
\hline
360 & 1000 & 600 & 320 & 680 & 102400 & 3686400S-4\\
400 & 900 & 600 & 250 & 650 & 62500 & 2250000S-4\\
450 & 800 & 600 & 175 & 625 & 30625 & 1102500S-4\\
480 & 750 & 600 & 135 & 615 & 18225 & 656100S-4\\
500 & 720 & 600 & 110 & 610 & 12100 & 435600S-4\\
\hline
\end{tabular}
\end{center}
\begin{center}
\begin{tabular}{ccccccc}
\hline
$x$ & $y$ & $b=12\cdot 80$  & $a$ & $d$ & $a^{2}$ & $a^{2}/b^{2}$  \\
\hline
2 & 460800 & 960 & 230399 & 230401 & 53083699201 & 2687362272050625S-6\\
4 & 230400 & 960 & 115198 & 115202 & 13270579204 & 671823072202500S-6\\
6 & 153600 & 960 & 76797 & 76803 & 5897779209 & 298575072455625S-6\\
8 & 115200 & 960 & 57596 & 57604 & 3317299216 & 167938272810000S-6\\
10 & 92160 & 960 & 46075 & 46085 & 2122905625 & 107472097265625S-6\\
12 & 76800 & 960 & 38394 & 38406 & 1474099236 & 74626273822500S-6\\
16 & 57600 & 960 & 28792 & 28808 & 828979264 & 41967075240000S-6\\
18 & 51200 & 960 & 25591 & 25609 & 654899281 & 33154276100625S-6\\
20 & 46080 & 960 & 23030 & 23050 & 530380900 & 26850533062500S-6\\
24 & 38400 & 960 & 19188 & 19212 & 368179344 & 18639079290000S-6\\
30 & 30720 & 960 & 15345 & 15375 & 235469025 & 11920619390625S-6\\
32 & 28800 & 960 & 14384 & 14416 & 206899456 & 10474284960000S-6\\
36 & 25600 & 960 & 12782 & 12818 & 163379524 & 8271088402500S-6\\
40 & 23040 & 960 & 11500 & 11540 & 132250000 & 6695156250000S-6\\
48 & 19200 & 960 & 9576 & 9624 & 91699776 & 4642301160000S-6\\
50 & 18432 & 960 & 9191 & 9241 & 84474481 & 4276520600625S-6\\
60 & 15360 & 960 & 7650 & 7710 & 58522500 & 2962701562500S-6\\
64 & 14400 & 960 & 7168 & 7232 & 51380224 & 2601123840000S-6\\
72 & 12800 & 960 & 6364 & 6436 & 40500496 & 2050337610000S-6\\
80 & 11520 & 960 & 5720 & 5800 & 32718400 & 1656369000000S-6\\
90 & 10240 & 960 & 5075 & 5165 & 25755625 & 1303878515625S-6\\
96 & 9600 & 960 & 4752 & 4848 & 22581504 & 1143188640000S-6\\
100 & 9216 & 960 & 4558 & 4658 & 20775364 & 1051752802500S-6\\
120 & 7680 & 960 & 3780 & 3900 & 14288400 & 723350250000S-6\\
128 & 7200 & 960 & 3536 & 3664 & 12503296 & 632979360000S-6\\
144 & 6400 & 960 & 3128 & 3272 & 9784384 & 495334440000S-6\\
150 & 6144 & 960 & 2997 & 3147 & 8982009 & 454714205625S-6\\
160 & 5760 & 960 & 2800 & 2960 & 7840000 & 396900000000S-6\\
180 & 5120 & 960 & 2470 & 2650 & 6100900 & 308858062500S-6\\
192 & 4800 & 960 & 2304 & 2496 & 5308416 & 268738560000S-6\\
200 & 4608 & 960 & 2204 & 2404 & 4857616 & 245916810000S-6\\
240 & 3840 & 960 & 1800 & 2040 & 3240000 & 164025000000S-6\\
256 & 3600 & 960 & 1672 & 1928 & 2795584 & 141526440000S-6\\
288 & 3200 & 960 & 1456 & 1744 & 2119936 & 107321760000S-6\\
300 & 3072 & 960 & 1386 & 1686 & 1920996 & 97250422500S-6\\
320 & 2880 & 960 & 1280 & 1600 & 1638400 & 82944000000S-6\\
360 & 2560 & 960 & 1100 & 1460 & 1210000 & 61256250000S-6\\
384 & 2400 & 960 & 1008 & 1392 & 1016064 & 51438240000S-6\\
\hline
400 & 2304 & 960 & 952 & 1352 & 906304 & 45881640000S-6\\
450 & 2048 & 960 & 799 & 1249 & 638401 & 32319050625S-6\\
480 & 1920 & 960 & 720 & 1200 & 518400 & 26244000000S-6\\
512 & 1800 & 960 & 644 & 1156 & 414736 & 20996010000S-6\\
\hline
576 & 1600 & 960 & 512 & 1088 & 262144 & 13271040000S-6\\
600 & 1536 & 960 & 468 & 1068 & 219024 & 11088090000S-6\\
640 & 1440 & 960 & 400 & 1040 & 160000 & 8100000000S-6\\
720 & 1280 & 960 & 280 & 1000 & 78400 & 3969000000S-6\\
768 & 1200 & 960 & 216 & 984 & 46656 & 2361960000S-6\\
800 & 1152 & 960 & 176 & 976 & 30976 & 1568160000S-6\\
900 & 1024 & 960 & 62 & 962 & 3844 & 194602500S-6\\
\hline
\end{tabular}
\end{center}

\begin{center}
\begin{tabular}{ccccccc}
\hline
$x$ & $y$ & $b=12\cdot 200$ & $a$ & $d$ & $a^{2}$ & $a^{2}/b^{2}$  \\
\hline
2 & 2880000 & 2400 & 1439999 & 1440001 & 2073597120001 & 16796136672008100S-6\\
4 & 1440000 & 2400 & 719998 & 720002 & 518397120004 & 4199016672032400S-6\\
6 & 960000 & 2400 & 479997 & 480003 & 230397120009 & 1866216672072900S-6\\
8 & 720000 & 2400 & 359996 & 360004 & 129597120016 & 1049736672129600S-6\\
10 & 576000 & 2400 & 287995 & 288005 & 82941120025 & 671823072202500S-6\\
12 & 480000 & 2400 & 239994 & 240006 & 57597120036 & 466536672291600S-6\\
16 & 360000 & 2400 & 179992 & 180008 & 32397120064 & 262416672518400S-6\\
18 & 320000 & 2400 & 159991 & 160009 & 25597120081 & 207336672656100S-6\\
20 & 288000 & 2400 & 143990 & 144010 & 20733120100 & 167938272810000S-6\\
24 & 240000 & 2400 & 119988 & 120012 & 14397120144 & 116616673166400S-6\\
30 & 192000 & 2400 & 95985 & 96015 & 9213120225 & 74626273822500S-6\\
32 & 180000 & 2400 & 89984 & 90016 & 8097120256 & 65586674073600S-6\\
36 & 160000 & 2400 & 79982 & 80018 & 6397120324 & 51816674624400S-6\\
40 & 144000 & 2400 & 71980 & 72020 & 5181120400 & 41967075240000S-6\\
48 & 120000 & 2400 & 59976 & 60024 & 3597120576 & 29136676665600S-6\\
50 & 115200 & 2400 & 57575 & 57625 & 3314880625 & 26850533062500S-6\\
60 & 96000 & 2400 & 47970 & 48030 & 2301120900 & 18639079290000S-6\\
64 & 90000 & 2400 & 44968 & 45032 & 2022121024 & 16379180294400S-6\\
72 & 80000 & 2400 & 39964 & 40036 & 1597121296 & 12936682497600S-6\\
80 & 72000 & 2400 & 35960 & 36040 & 1293121600 & 10474284960000S-6\\
90 & 64000 & 2400 & 31955 & 32045 & 1021122025 & 8271088402500S-6\\
96 & 60000 & 2400 & 29952 & 30048 & 897122304 & 7266690662400S-6\\
100 & 57600 & 2400 & 28750 & 28850 & 826562500 & 6695156250000S-6\\
120 & 48000 & 2400 & 23940 & 24060 & 573123600 & 4642301160000S-6\\
128 & 45000 & 2400 & 22436 & 22564 & 503374096 & 4077330177600S-6\\
144 & 40000 & 2400 & 19928 & 20072 & 397125184 & 3216713990400S-6\\
150 & 38400 & 2400 & 19125 & 19275 & 365765625 & 2962701562500S-6\\
160 & 36000 & 2400 & 17920 & 18080 & 321126400 & 2601123840000S-6\\
180 & 32000 & 2400 & 15910 & 16090 & 253128100 & 2050337610000S-6\\
192 & 30000 & 2400 & 14904 & 15096 & 222129216 & 1799246649600S-6\\
200 & 28800 & 2400 & 14300 & 14500 & 204490000 & 1656369000000S-6\\
240 & 24000 & 2400 & 11880 & 12120 & 141134400 & 1143188640000S-6\\
250 & 23040 & 2400 & 11395 & 11645 & 129846025 & 1051752802500S-6\\
256 & 22500 & 2400 & 11122 & 11378 & 123698884 & 1001960960400S-6\\
288 & 20000 & 2400 & 9856 & 10144 & 97140736 & 786839961600S-6\\
300 & 19200 & 2400 & 9450 & 9750 & 89302500 & 723350250000S-6\\
320 & 18000 & 2400 & 8840 & 9160 & 78145600 & 632979360000S-6\\
360 & 16000 & 2400 & 7820 & 8180 & 61152400 & 495334440000S-6\\
384 & 15000 & 2400 & 7308 & 7692 & 53406864 & 432595598400S-6\\
400 & 14400 & 2400 & 7000 & 7400 & 49000000 & 396900000000S-6\\
450 & 12800 & 2400 & 6175 & 6625 & 38130625 & 308858062500S-6\\
480 & 12000 & 2400 & 5760 & 6240 & 33177600 & 268738560000S-6\\
500 & 11520 & 2400 & 5510 & 6010 & 30360100 & 245916810000S-6\\
512 & 11250 & 2400 & 5369 & 5881 & 28826161 & 233491904100S-6\\
576 & 10000 & 2400 & 4712 & 5288 & 22202944 & 179843846400S-6\\
600 & 9600 & 2400 & 4500 & 5100 & 20250000 & 164025000000S-6\\
\hline
\end{tabular}
\end{center}
\begin{center}
\begin{tabular}{ccccccc}
\hline
$x$ & $y$ & $b=12\cdot 200$ & $a$ & $d$ & $a^{2}$ & $a^{2}/b^{2}$  \\
\hline
640 & 9000 & 2400 & 4180 & 4820 & 17472400 & 141526440000S-6\\
720 & 8000 & 2400 & 3640 & 4360 & 13249600 & 107321760000S-6\\
750 & 7680 & 2400 & 3465 & 4215 & 12006225 & 97250422500S-6\\
768 & 7500 & 2400 & 3366 & 4134 & 11329956 & 91772643600S-6\\
800 & 7200 & 2400 & 3200 & 4000 & 10240000 & 82944000000S-6\\
900 & 6400 & 2400 & 2750 & 3650 & 7562500 & 61256250000S-6\\
960 & 6000 & 2400 & 2520 & 3480 & 6350400 & 51438240000S-6\\
\hline
1000 & 5760 & 2400 & 2380 & 3380 & 5664400 & 45881640000S-6\\
1152 & 5000 & 2400 & 1924 & 3076 & 3701776 & 29984385600S-6\\
1200 & 4800 & 2400 & 1800 & 3000 & 3240000 & 26244000000S-6\\
1250 & 4608 & 2400 & 1679 & 2929 & 2819041 & 22834232100S-6\\
1280 & 4500 & 2400 & 1610 & 2890 & 2592100 & 20996010000S-6\\
\hline
1440 & 4000 & 2400 & 1280 & 2720 & 1638400 & 13271040000S-6\\
1500 & 3840 & 2400 & 1170 & 2670 & 1368900 & 11088090000S-6\\
1536 & 3750 & 2400 & 1107 & 2643 & 1225449 & 9926136900S-6\\
1600 & 3600 & 2400 & 1000 & 2600 & 1000000 & 8100000000S-6\\
1800 & 3200 & 2400 & 700 & 2500 & 490000 & 3969000000S-6\\
1920 & 3000 & 2400 & 540 & 2460 & 291600 & 2361960000S-6\\
2000 & 2880 & 2400 & 440 & 2440 & 193600 & 1568160000S-6\\
2250 & 2560 & 2400 & 155 & 2405 & 24025 & 194602500S-6\\
2304 & 2500 & 2400 & 98 & 2402 & 9604 & 77792400S-6\\
\hline
\end{tabular}
\end{center}
\begin{center}
\begin{tabular}{ccccccc}
\hline
$x$ & $y$ & $b=12\cdot 225$ & $a$ & $d$ & $a^{2}$ & $a^{2}/b^{2}$  \\
\hline
2 & 3645000 & 2700 & 1822499 & 1822501 & 3321502605001 & 21257616672006400S-6\\
4 & 1822500 & 2700 & 911248 & 911252 & 830372917504 & 5314386672025600S-6\\
6 & 1215000 & 2700 & 607497 & 607503 & 369052605009 & 2361936672057600S-6\\
8 & 911250 & 2700 & 455621 & 455629 & 207590495641 & 1328579172102400S-6\\
10 & 729000 & 2700 & 364495 & 364505 & 132856605025 & 850282272160000S-6\\
12 & 607500 & 2700 & 303744 & 303756 & 92260417536 & 590466672230400S-6\\
18 & 405000 & 2700 & 202491 & 202509 & 41002605081 & 262416672518400S-6\\
20 & 364500 & 2700 & 182240 & 182260 & 33211417600 & 212553072640000S-6\\
24 & 303750 & 2700 & 151863 & 151887 & 23062370769 & 147599172921600S-6\\
30 & 243000 & 2700 & 121485 & 121515 & 14758605225 & 94455073440000S-6\\
36 & 202500 & 2700 & 101232 & 101268 & 10247917824 & 65586674073600S-6\\
40 & 182250 & 2700 & 91105 & 91145 & 8300121025 & 53120774560000S-6\\
50 & 145800 & 2700 & 72875 & 72925 & 5310765625 & 33988900000000S-6\\
54 & 135000 & 2700 & 67473 & 67527 & 4552605729 & 29136676665600S-6\\
60 & 121500 & 2700 & 60720 & 60780 & 3686918400 & 23596277760000S-6\\
72 & 101250 & 2700 & 50589 & 50661 & 2559246921 & 16379180294400S-6\\
90 & 81000 & 2700 & 40455 & 40545 & 1636607025 & 10474284960000S-6\\
100 & 72900 & 2700 & 36400 & 36500 & 1324960000 & 8479744000000S-6\\
108 & 67500 & 2700 & 33696 & 33804 & 1135420416 & 7266690662400S-6\\
120 & 60750 & 2700 & 30315 & 30435 & 918999225 & 5881595040000S-6\\
150 & 48600 & 2700 & 24225 & 24375 & 586850625 & 3755844000000S-6\\
\hline
\end{tabular}
\end{center}
\begin{center}
\begin{tabular}{ccccccc}
\hline
$x$ & $y$ & $b=12\cdot 225$ & $a$ & $d$ & $a^{2}$ & $a^{2}/b^{2}$  \\
\hline
162 & 45000 & 2700 & 22419 & 22581 & 502611561 & 3216713990400S-6\\
180 & 40500 & 2700 & 20160 & 20340 & 406425600 & 2601123840000S-6\\
200 & 36450 & 2700 & 18125 & 18325 & 328515625 & 2102500000000S-6\\
216 & 33750 & 2700 & 16767 & 16983 & 281132289 & 1799246649600S-6\\
250 & 29160 & 2700 & 14455 & 14705 & 208947025 & 1337260960000S-6\\
270 & 27000 & 2700 & 13365 & 13635 & 178623225 & 1143188640000S-6\\
300 & 24300 & 2700 & 12000 & 12300 & 144000000 & 921600000000S-6\\
324 & 22500 & 2700 & 11088 & 11412 & 122943744 & 786839961600S-6\\
360 & 20250 & 2700 & 9945 & 10305 & 98903025 & 632979360000S-6\\
450 & 16200 & 2700 & 7875 & 8325 & 62015625 & 396900000000S-6\\
486 & 15000 & 2700 & 7257 & 7743 & 52664049 & 337049913600S-6\\
500 & 14580 & 2700 & 7040 & 7540 & 49561600 & 317194240000S-6\\
540 & 13500 & 2700 & 6480 & 7020 & 41990400 & 268738560000S-6\\
600 & 12150 & 2700 & 5775 & 6375 & 33350625 & 213444000000S-6\\
648 & 11250 & 2700 & 5301 & 5949 & 28100601 & 179843846400S-6\\
750 & 9720 & 2700 & 4485 & 5235 & 20115225 & 128737440000S-6\\
810 & 9000 & 2700 & 4095 & 4905 & 16769025 & 107321760000S-6\\
900 & 8100 & 2700 & 3600 & 4500 & 12960000 & 82944000000S-6\\
972 & 7500 & 2700 & 3264 & 4236 & 10653696 & 68183654400S-6\\
1000 & 7290 & 2700 & 3145 & 4145 & 9891025 & 63302560000S-6\\
1080 & 6750 & 2700 & 2835 & 3915 & 8037225 & 51438240000S-6\\
\hline
1250 & 5832 & 2700 & 2291 & 3541 & 5248681 & 33591558400S-6\\
1350 & 5400 & 2700 & 2025 & 3375 & 4100625 & 26244000000S-6\\
1458 & 5000 & 2700 & 1771 & 3229 & 3136441 & 20073222400S-6\\
1500 & 4860 & 2700 & 1680 & 3180 & 2822400 & 18063360000S-6\\
\hline
1620 & 4500 & 2700 & 1440 & 3060 & 2073600 & 13271040000S-6\\
1800 & 4050 & 2700 & 1125 & 2925 & 1265625 & 8100000000S-6\\
1944 & 3750 & 2700 & 903 & 2847 & 815409 & 5218617600S-6\\
2250 & 3240 & 2700 & 495 & 2745 & 245025 & 1568160000S-6\\
2430 & 3000 & 2700 & 285 & 2715 & 81225 & 519840000S-6\\
2500 & 2916 & 2700 & 208 & 2708 & 43264 & 276889600S-6\\
\hline
\end{tabular}
\end{center}
\begin{center}
\begin{tabular}{ccccccc}
\hline
$x$ & $y$ & $a$ & $b=12\cdot 288$ & $d$ & $b^{2}$ & $b^{2}/a^{2}$ \\
\hline
2 & 5971968 & 3456 & 2985983 & 2985985 & 8916094476289 & 34273467166854916S-6\\
4 & 2985984 & 3456 & 1492990 & 1492994 & 2229019140100 & 8568349574544400S-6\\
6 & 1990656 & 3456 & 995325 & 995331 & 990671855625 & 3808142613022500S-6\\
8 & 1492992 & 3456 & 746492 & 746500 & 557250306064 & 2142070176510016S-6\\
12 & 995328 & 3456 & 497658 & 497670 & 247663484964 & 952018436201616S-6\\
16 & 746496 & 3456 & 373240 & 373256 & 139308097600 & 535500327174400S-6\\
18 & 663552 & 3456 & 331767 & 331785 & 110069342289 & 423106551758916S-6\\
24 & 497664 & 3456 & 248820 & 248844 & 61911392400 & 237987392385600S-6\\
32 & 373248 & 3456 & 186608 & 186640 & 34822545664 & 133857865532416S-6\\
36 & 331776 & 3456 & 165870 & 165906 & 27512856900 & 105759421923600S-6\\
48 & 248832 & 3456 & 124392 & 124440 & 15473369664 & 59479632988416S-6\\
54 & 221184 & 3456 & 110565 & 110619 & 12224619225 & 46991436300900S-6\\
64 & 186624 & 3456 & 93280 & 93344 & 8701158400 & 33447252889600S-6\\
72 & 165888 & 3456 & 82908 & 82980 & 6873736464 & 26422642967616S-6\\
96 & 124416 & 3456 & 62160 & 62256 & 3863865600 & 14852699366400S-6\\
108 & 110592 & 3456 & 55242 & 55350 & 3051678564 & 11730652400016S-6\\
128 & 93312 & 3456 & 46592 & 46720 & 2170814464 & 8344610799616S-6\\
144 & 82944 & 3456 & 41400 & 41544 & 1713960000 & 6588462240000S-6\\
162 & 73728 & 3456 & 36783 & 36945 & 1352989089 & 5200890058116S-6\\
192 & 62208 & 3456 & 31008 & 31200 & 961496064 & 3695990870016S-6\\
216 & 55296 & 3456 & 27540 & 27756 & 758451600 & 2915487950400S-6\\
256 & 46656 & 3456 & 23200 & 23456 & 538240000 & 2068994560000S-6\\
288 & 41472 & 3456 & 20592 & 20880 & 424030464 & 1629973103616S-6\\
324 & 36864 & 3456 & 18270 & 18594 & 333792900 & 1283099907600S-6\\
384 & 31104 & 3456 & 15360 & 15744 & 235929600 & 906913382400S-6\\
432 & 27648 & 3456 & 13608 & 14040 & 185177664 & 711822940416S-6\\
486 & 24576 & 3456 & 12045 & 12531 & 145082025 & 557695304100S-6\\
512 & 23328 & 3456 & 11408 & 11920 & 130142464 & 500267631616S-6\\
576 & 20736 & 3456 & 10080 & 10656 & 101606400 & 390575001600S-6\\
648 & 18432 & 3456 & 8892 & 9540 & 79067664 & 303936100416S-6\\
768 & 15552 & 3456 & 7392 & 8160 & 54641664 & 210042556416S-6\\
864 & 13824 & 3456 & 6480 & 7344 & 41990400 & 161411097600S-6\\
972 & 12288 & 3456 & 5658 & 6630 & 32012964 & 123057833616S-6\\
1024 & 11664 & 3456 & 5320 & 6344 & 28302400 & 108794425600S-6\\
1152 & 10368 & 3456 & 4608 & 5760 & 21233664 & 81622204416S-6\\
1296 & 9216 & 3456 & 3960 & 5256 & 15681600 & 60280070400S-6\\
\hline
1458 & 8192 & 3456 & 3367 & 4825 & 11336689 & 43578232516S-6\\
1536 & 7776 & 3456 & 3120 & 4656 & 9734400 & 37419033600S-6\\
1728 & 6912 & 3456 & 2592 & 4320 & 6718464 & 25825775616S-6\\
1944 & 6144 & 3456 & 2100 & 4044 & 4410000 & 16952040000S-6\\
\hline
2048 & 5832 & 3456 & 1892 & 3940 & 3579664 & 13760228416S-6\\
2304 & 5184 & 3456 & 1440 & 3744 & 2073600 & 7970918400S-6\\
2592 & 4608 & 3456 & 1008 & 3600 & 1016064 & 3905750016S-6\\
2916 & 4096 & 3456 & 590 & 3506 & 348100 & 1338096400S-6\\
3072 & 3888 & 3456 & 408 & 3480 & 166464 & 639887616S-6\\
\hline
\end{tabular}
\end{center}

\begin{center}
\begin{tabular}{ccccccc}
\hline
$x$ & $y$ & $b=12\cdot 400$ & $a$ & $d$ & $b^{2}$ & $b^{2}/a^{2}$ \\
\hline
2 & 11520000 & 4800 & 5759999 & 5760001 & 33177588480001 & 67184616672002025S-6\\
4 & 5760000 & 4800 & 2879998 & 2880002 & 8294388480004 & 16796136672008100S-6\\
6 & 3840000 & 4800 & 1919997 & 1920003 & 3686388480009 & 7464936672018225S-6\\
8 & 2880000 & 4800 & 1439996 & 1440004 & 2073588480016 & 4199016672032400S-6\\
10 & 2304000 & 4800 & 1151995 & 1152005 & 1327092480025 & 2687362272050625S-6\\
12 & 1920000 & 4800 & 959994 & 960006 & 921588480036 & 1866216672072900S-6\\
16 & 1440000 & 4800 & 719992 & 720008 & 518388480064 & 1049736672129600S-6\\
18 & 1280000 & 4800 & 639991 & 640009 & 409588480081 & 829416672164025S-6\\
20 & 1152000 & 4800 & 575990 & 576010 & 331764480100 & 671823072202500S-6\\
24 & 960000 & 4800 & 479988 & 480012 & 230388480144 & 466536672291600S-6\\
30 & 768000 & 4800 & 383985 & 384015 & 147444480225 & 298575072455625S-6\\
32 & 720000 & 4800 & 359984 & 360016 & 129588480256 & 262416672518400S-6\\
36 & 640000 & 4800 & 319982 & 320018 & 102388480324 & 207336672656100S-6\\
40 & 576000 & 4800 & 287980 & 288020 & 82932480400 & 167938272810000S-6\\
48 & 480000 & 4800 & 239976 & 240024 & 57588480576 & 116616673166400S-6\\
50 & 460800 & 4800 & 230375 & 230425 & 53072640625 & 107472097265625S-6\\
60 & 384000 & 4800 & 191970 & 192030 & 36852480900 & 74626273822500S-6\\
64 & 360000 & 4800 & 179968 & 180032 & 32388481024 & 65586674073600S-6\\
72 & 320000 & 4800 & 159964 & 160036 & 25588481296 & 51816674624400S-6\\
80 & 288000 & 4800 & 143960 & 144040 & 20724481600 & 41967075240000S-6\\
90 & 256000 & 4800 & 127955 & 128045 & 16372482025 & 33154276100625S-6\\
96 & 240000 & 4800 & 119952 & 120048 & 14388482304 & 29136676665600S-6\\
100 & 230400 & 4800 & 115150 & 115250 & 13259522500 & 26850533062500S-6\\
120 & 192000 & 4800 & 95940 & 96060 & 9204483600 & 18639079290000S-6\\
128 & 180000 & 4800 & 89936 & 90064 & 8088484096 & 16379180294400S-6\\
144 & 160000 & 4800 & 79928 & 80072 & 6388485184 & 12936682497600S-6\\
150 & 153600 & 4800 & 76725 & 76875 & 5886725625 & 11920619390625S-6\\
160 & 144000 & 4800 & 71920 & 72080 & 5172486400 & 10474284960000S-6\\
180 & 128000 & 4800 & 63910 & 64090 & 4084488100 & 8271088402500S-6\\
192 & 120000 & 4800 & 59904 & 60096 & 3588489216 & 7266690662400S-6\\
200 & 115200 & 4800 & 57500 & 57700 & 3306250000 & 6695156250000S-6\\
240 & 96000 & 4800 & 47880 & 48120 & 2292494400 & 4642301160000S-6\\
250 & 92160 & 4800 & 45955 & 46205 & 2111862025 & 4276520600625S-6\\
256 & 90000 & 4800 & 44872 & 45128 & 2013496384 & 4077330177600S-6\\
288 & 80000 & 4800 & 39856 & 40144 & 1588500736 & 3216713990400S-6\\
300 & 76800 & 4800 & 38250 & 38550 & 1463062500 & 2962701562500S-6\\
320 & 72000 & 4800 & 35840 & 36160 & 1284505600 & 2601123840000S-6\\
360 & 64000 & 4800 & 31820 & 32180 & 1012512400 & 2050337610000S-6\\
384 & 60000 & 4800 & 29808 & 30192 & 888516864 & 1799246649600S-6\\
400 & 57600 & 4800 & 28600 & 29000 & 817960000 & 1656369000000S-6\\
450 & 51200 & 4800 & 25375 & 25825 & 643890625 & 1303878515625S-6\\
480 & 48000 & 4800 & 23760 & 24240 & 564537600 & 1143188640000S-6\\
500 & 46080 & 4800 & 22790 & 23290 & 519384100 & 1051752802500S-6\\
512 & 45000 & 4800 & 22244 & 22756 & 494795536 & 1001960960400S-6\\
576 & 40000 & 4800 & 19712 & 20288 & 388562944 & 786839961600S-6\\
600 & 38400 & 4800 & 18900 & 19500 & 357210000 & 723350250000S-6\\
\hline
\end{tabular}
\end{center}
\begin{center}
\begin{tabular}{ccccccc}
\hline
$x$ & $y$ & $b=12\cdot 400$ & $a$ & $d$ & $b^{2}$ & $b^{2}/a^{2}$ \\
\hline
640 & 36000 & 4800 & 17680 & 18320 & 312582400 & 632979360000S-6\\
720 & 32000 & 4800 & 15640 & 16360 & 244609600 & 495334440000S-6\\
750 & 30720 & 4800 & 14985 & 15735 & 224550225 & 454714205625S-6\\
768 & 30000 & 4800 & 14616 & 15384 & 213627456 & 432595598400S-6\\
800 & 28800 & 4800 & 14000 & 14800 & 196000000 & 396900000000S-6\\
900 & 25600 & 4800 & 12350 & 13250 & 152522500 & 308858062500S-6\\
960 & 24000 & 4800 & 11520 & 12480 & 132710400 & 268738560000S-6\\
1000 & 23040 & 4800 & 11020 & 12020 & 121440400 & 245916810000S-6\\
1024 & 22500 & 4800 & 10738 & 11762 & 115304644 & 233491904100S-6\\
1152 & 20000 & 4800 & 9424 & 10576 & 88811776 & 179843846400S-6\\
1200 & 19200 & 4800 & 9000 & 10200 & 81000000 & 164025000000S-6\\
1250 & 18432 & 4800 & 8591 & 9841 & 73805281 & 149455694025S-6\\
1280 & 18000 & 4800 & 8360 & 9640 & 69889600 & 141526440000S-6\\
1440 & 16000 & 4800 & 7280 & 8720 & 52998400 & 107321760000S-6\\
1500 & 15360 & 4800 & 6930 & 8430 & 48024900 & 97250422500S-6\\
1536 & 15000 & 4800 & 6732 & 8268 & 45319824 & 91772643600S-6\\
1600 & 14400 & 4800 & 6400 & 8000 & 40960000 & 82944000000S-6\\
1800 & 12800 & 4800 & 5500 & 7300 & 30250000 & 61256250000S-6\\
1920 & 12000 & 4800 & 5040 & 6960 & 25401600 & 51438240000S-6\\
\hline
2000 & 11520 & 4800 & 4760 & 6760 & 22657600 & 45881640000S-6\\
2048 & 11250 & 4800 & 4601 & 6649 & 21169201 & 42867632025S-6\\
2250 & 10240 & 4800 & 3995 & 6245 & 15960025 & 32319050625S-6\\
2304 & 10000 & 4800 & 3848 & 6152 & 14807104 & 29984385600S-6\\
2400 & 9600 & 4800 & 3600 & 6000 & 12960000 & 26244000000S-6\\
2500 & 9216 & 4800 & 3358 & 5858 & 11276164 & 22834232100S-6\\
2560 & 9000 & 4800 & 3220 & 5780 & 10368400 & 20996010000S-6\\
\hline
2880 & 8000 & 4800 & 2560 & 5440 & 6553600 & 13271040000S-6\\
3000 & 7680 & 4800 & 2340 & 5340 & 5475600 & 11088090000S-6\\
3072 & 7500 & 4800 & 2214 & 5286 & 4901796 & 9926136900S-6\\
3200 & 7200 & 4800 & 2000 & 5200 & 4000000 & 8100000000S-6\\
3600 & 6400 & 4800 & 1400 & 5000 & 1960000 & 3969000000S-6\\
3750 & 6144 & 4800 & 1197 & 4947 & 1432809 & 2901438225S-6\\
3840 & 6000 & 4800 & 1080 & 4920 & 1166400 & 2361960000S-6\\
4000 & 5760 & 4800 & 880 & 4880 & 774400 & 1568160000S-6\\
4500 & 5120 & 4800 & 310 & 4810 & 96100 & 194602500S-6\\
4608 & 5000 & 4800 & 196 & 4804 & 38416 & 77792400S-6\\
\hline
\end{tabular}
\end{center}
\begin{center}
\begin{tabular}{ccccccc}
\hline
$x$ & $y$ & $b=12\cdot 540$ & $a$ & $d$ & $b^{2}$ & $b^{2}/a^{2}$ \\
\hline
2 & 20995200 & 6480 & 10497599 & 10497601 & 110199584764801 & 16523225363884312832S-8\\
4 & 10497600 & 6480 & 5248798 & 5248802 & 27549880444804 & 17965801410668241920S-8\\
6 & 6998400 & 6480 & 3499197 & 3499203 & 12244379644809 & 12084030431816896768S-8\\
8 & 5248800 & 6480 & 2624396 & 2624404 & 6887454364816 & 9103073385554448384S-8\\
10 & 4199040 & 6480 & 2099515 & 2099525 & 4407963235225 & 17631852940900000000S-8\\
12 & 3499200 & 6480 & 1749594 & 1749606 & 3061079164836 & 12244316659344000000S-8\\
16 & 2624400 & 6480 & 1312192 & 1312208 & 1721847844864 & 6887391379456000000S-8\\
18 & 2332800 & 6480 & 1166391 & 1166409 & 1360467964881 & 5441871859524000000S-8\\
20 & 2099520 & 6480 & 1049750 & 1049770 & 1101975062500 & 4407900250000000000S-8\\
24 & 1749600 & 6480 & 874788 & 874812 & 765254044944 & 3061016179776000000S-8\\
30 & 1399680 & 6480 & 699825 & 699855 & 489755030625 & 1959020122500000000S-8\\
32 & 1312200 & 6480 & 656084 & 656116 & 430446215056 & 1721784860224000000S-8\\
36 & 1166400 & 6480 & 583182 & 583218 & 340101245124 & 1360404980496000000S-8\\
40 & 1049760 & 6480 & 524860 & 524900 & 275478019600 & 1101912078400000000S-8\\
48 & 874800 & 6480 & 437376 & 437424 & 191297765376 & 765191061504000000S-8\\
50 & 839808 & 6480 & 419879 & 419929 & 176298374641 & 705193498564000000S-8\\
54 & 777600 & 6480 & 388773 & 388827 & 151144445529 & 604577782116000000S-8\\
60 & 699840 & 6480 & 349890 & 349950 & 122423012100 & 489692048400000000S-8\\
64 & 656100 & 6480 & 328018 & 328082 & 107595808324 & 430383233296000000S-8\\
72 & 583200 & 6480 & 291564 & 291636 & 85009566096 & 340038264384000000S-8\\
80 & 524880 & 6480 & 262400 & 262480 & 68853760000 & 275415040000000000S-8\\
90 & 466560 & 6480 & 233235 & 233325 & 54398565225 & 217594260900000000S-8\\
96 & 437400 & 6480 & 218652 & 218748 & 47808697104 & 191234788416000000S-8\\
100 & 419904 & 6480 & 209902 & 210002 & 44058849604 & 176235398416000000S-8\\
108 & 388800 & 6480 & 194346 & 194454 & 37770367716 & 151081470864000000S-8\\
120 & 349920 & 6480 & 174900 & 175020 & 30590010000 & 122360040000000000S-8\\
128 & 328050 & 6480 & 163961 & 164089 & 26883209521 & 107532838084000000S-8\\
144 & 291600 & 6480 & 145728 & 145872 & 21236649984 & 84946599936000000S-8\\
150 & 279936 & 6480 & 139893 & 140043 & 19570051449 & 78280205796000000S-8\\
160 & 262440 & 6480 & 131140 & 131300 & 17197699600 & 68790798400000000S-8\\
162 & 259200 & 6480 & 129519 & 129681 & 16775171361 & 67100685444000000S-8\\
180 & 233280 & 6480 & 116550 & 116730 & 13583902500 & 54335610000000000S-8\\
192 & 218700 & 6480 & 109254 & 109446 & 11936436516 & 47745746064000000S-8\\
200 & 209952 & 6480 & 104876 & 105076 & 10998975376 & 43995901504000000S-8\\
216 & 194400 & 6480 & 97092 & 97308 & 9426856464 & 37707425856000000S-8\\
240 & 174960 & 6480 & 87360 & 87600 & 7631769600 & 30527078400000000S-8\\
270 & 155520 & 6480 & 77625 & 77895 & 6025640625 & 24102562500000000S-8\\
288 & 145800 & 6480 & 72756 & 73044 & 5293435536 & 21173742144000000S-8\\
300 & 139968 & 6480 & 69834 & 70134 & 4876787556 & 19507150224000000S-8\\
320 & 131220 & 6480 & 65450 & 65770 & 4283702500 & 17134810000000000S-8\\
324 & 129600 & 6480 & 64638 & 64962 & 4178071044 & 16712284176000000S-8\\
360 & 116640 & 6480 & 58140 & 58500 & 3380259600 & 13521038400000000S-8\\
384 & 109350 & 6480 & 54483 & 54867 & 2968397289 & 11873589156000000S-8\\
400 & 104976 & 6480 & 52288 & 52688 & 2734034944 & 10936139776000000S-8\\
432 & 97200 & 6480 & 48384 & 48816 & 2341011456 & 9364045824000000S-8\\
450 & 93312 & 6480 & 46431 & 46881 & 2155837761 & 8623351044000000S-8\\
480 & 87480 & 6480 & 43500 & 43980 & 1892250000 & 7569000000000000S-8\\
\hline
\end{tabular}
\end{center}
\begin{center}
\begin{tabular}{ccccccc}
\hline
$x$ & $y$ & $b=12\cdot 540$ & $a$ & $d$ & $b^{2}$ & $b^{2}/a^{2}$ \\
\hline
486 & 86400 & 6480 & 42957 & 43443 & 1845303849 & 7381215396000000S-8\\
540 & 77760 & 6480 & 38610 & 39150 & 1490732100 & 5962928400000000S-8\\
576 & 72900 & 6480 & 36162 & 36738 & 1307690244 & 5230760976000000S-8\\
600 & 69984 & 6480 & 34692 & 35292 & 1203534864 & 4814139456000000S-8\\
640 & 65610 & 6480 & 32485 & 33125 & 1055275225 & 4221100900000000S-8\\
648 & 64800 & 6480 & 32076 & 32724 & 1028869776 & 4115479104000000S-8\\
720 & 58320 & 6480 & 28800 & 29520 & 829440000 & 3317760000000000S-8\\
800 & 52488 & 6480 & 25844 & 26644 & 667912336 & 2671649344000000S-8\\
810 & 51840 & 6480 & 25515 & 26325 & 651015225 & 2604060900000000S-8\\
864 & 48600 & 6480 & 23868 & 24732 & 569681424 & 2278725696000000S-8\\
900 & 46656 & 6480 & 22878 & 23778 & 523402884 & 2093611536000000S-8\\
960 & 43740 & 6480 & 21390 & 22350 & 457532100 & 1830128400000000S-8\\
972 & 43200 & 6480 & 21114 & 22086 & 445800996 & 1783203984000000S-8\\
1080 & 38880 & 6480 & 18900 & 19980 & 357210000 & 1428840000000000S-8\\
1152 & 36450 & 6480 & 17649 & 18801 & 311487201 & 1245948804000000S-8\\
1200 & 34992 & 6480 & 16896 & 18096 & 285474816 & 1141899264000000S-8\\
1296 & 32400 & 6480 & 15552 & 16848 & 241864704 & 967458816000000S-8\\
1350 & 31104 & 6480 & 14877 & 16227 & 221325129 & 885300516000000S-8\\
1440 & 29160 & 6480 & 13860 & 15300 & 192099600 & 768398400000000S-8\\
1458 & 28800 & 6480 & 13671 & 15129 & 186896241 & 747584964000000S-8\\
1600 & 26244 & 6480 & 12322 & 13922 & 151831684 & 607326736000000S-8\\
1620 & 25920 & 6480 & 12150 & 13770 & 147622500 & 590490000000000S-8\\
1728 & 24300 & 6480 & 11286 & 13014 & 127373796 & 509495184000000S-8\\
1800 & 23328 & 6480 & 10764 & 12564 & 115863696 & 463454784000000S-8\\
1920 & 21870 & 6480 & 9975 & 11895 & 99500625 & 398002500000000S-8\\
1944 & 21600 & 6480 & 9828 & 11772 & 96589584 & 386358336000000S-8\\
2160 & 19440 & 6480 & 8640 & 10800 & 74649600 & 298598400000000S-8\\
2400 & 17496 & 6480 & 7548 & 9948 & 56972304 & 227889216000000S-8\\
2430 & 17280 & 6480 & 7425 & 9855 & 55130625 & 220522500000000S-8\\
2592 & 16200 & 6480 & 6804 & 9396 & 46294416 & 185177664000000S-8\\
\hline
2700 & 15552 & 6480 & 6426 & 9126 & 41293476 & 165173904000000S-8\\
2880 & 14580 & 6480 & 5850 & 8730 & 34222500 & 136890000000000S-8\\
2916 & 14400 & 6480 & 5742 & 8658 & 32970564 & 131882256000000S-8\\
3200 & 13122 & 6480 & 4961 & 8161 & 24611521 & 98446084000000S-8\\
3240 & 12960 & 6480 & 4860 & 8100 & 23619600 & 94478400000000S-8\\
3456 & 12150 & 6480 & 4347 & 7803 & 18896409 & 75585636000000S-8\\
3600 & 11664 & 6480 & 4032 & 7632 & 16257024 & 65028096000000S-8\\
\hline
3888 & 10800 & 6480 & 3456 & 7344 & 11943936 & 47775744000000S-8\\
4050 & 10368 & 6480 & 3159 & 7209 & 9979281 & 39917124000000S-8\\
4320 & 9720 & 6480 & 2700 & 7020 & 7290000 & 29160000000000S-8\\
4374 & 9600 & 6480 & 2613 & 6987 & 6827769 & 27311076000000S-8\\
4800 & 8748 & 6480 & 1974 & 6774 & 3896676 & 15586704000000S-8\\
4860 & 8640 & 6480 & 1890 & 6750 & 3572100 & 14288400000000S-8\\
5184 & 8100 & 6480 & 1458 & 6642 & 2125764 & 8503056000000S-8\\
5400 & 7776 & 6480 & 1188 & 6588 & 1411344 & 5645376000000S-8\\
5760 & 7290 & 6480 & 765 & 6525 & 585225 & 2340900000000S-8\\
5832 & 7200 & 6480 & 684 & 6516 & 467856 & 1871424000000S-8\\
\hline
\end{tabular}
\end{center}

\begin{center}
\begin{tabular}{ccccccc}
\hline
$x$ & $y$ & $b = 12\cdot 1125$ & $a$ & $d$ & $b^{2}$ & $b^{2}/a^{2}$ \\
\hline
2 & 91125000 & 13500 & 45562499 & 45562501 & 2075941315125001 & 531440976672000256S-6\\
4 & 45562500 & 13500 & 22781248 & 22781252 & 518985260437504 & 132860226672001024S-6\\
6 & 30375000 & 13500 & 15187497 & 15187503 & 230660065125009 & 59048976672002304S-6\\
8 & 22781250 & 13500 & 11390621 & 11390629 & 129746246765641 & 33215039172004096S-6\\
10 & 18225000 & 13500 & 9112495 & 9112505 & 83037565125025 & 21257616672006400S-6\\
12 & 15187500 & 13500 & 7593744 & 7593756 & 57664947937536 & 14762226672009216S-6\\
18 & 10125000 & 13500 & 5062491 & 5062509 & 25628815125081 & 6560976672020736S-6\\
20 & 9112500 & 13500 & 4556240 & 4556260 & 20759322937600 & 5314386672025600S-6\\
24 & 7593750 & 13500 & 3796863 & 3796887 & 14416168640769 & 3690539172036864S-6\\
30 & 6075000 & 13500 & 3037485 & 3037515 & 9226315125225 & 2361936672057600S-6\\
36 & 5062500 & 13500 & 2531232 & 2531268 & 6407135437824 & 1640226672082944S-6\\
40 & 4556250 & 13500 & 2278105 & 2278145 & 5189762391025 & 1328579172102400S-6\\
50 & 3645000 & 13500 & 1822475 & 1822525 & 3321415125625 & 850282272160000S-6\\
54 & 3375000 & 13500 & 1687473 & 1687527 & 2847565125729 & 728976672186624S-6\\
60 & 3037500 & 13500 & 1518720 & 1518780 & 2306510438400 & 590466672230400S-6\\
72 & 2531250 & 13500 & 1265589 & 1265661 & 1601715516921 & 410039172331776S-6\\
90 & 2025000 & 13500 & 1012455 & 1012545 & 1025065127025 & 262416672518400S-6\\
100 & 1822500 & 13500 & 911200 & 911300 & 830285440000 & 212553072640000S-6\\
108 & 1687500 & 13500 & 843696 & 843804 & 711822940416 & 182226672746496S-6\\
120 & 1518750 & 13500 & 759315 & 759435 & 576559269225 & 147599172921600S-6\\
150 & 1215000 & 13500 & 607425 & 607575 & 368965130625 & 94455073440000S-6\\
162 & 1125000 & 13500 & 562419 & 562581 & 316315131561 & 80976673679616S-6\\
180 & 1012500 & 13500 & 506160 & 506340 & 256197945600 & 65586674073600S-6\\
200 & 911250 & 13500 & 455525 & 455725 & 207503025625 & 53120774560000S-6\\
216 & 843750 & 13500 & 421767 & 421983 & 177887402289 & 45539174985984S-6\\
250 & 729000 & 13500 & 364375 & 364625 & 132769140625 & 33988900000000S-6\\
270 & 675000 & 13500 & 337365 & 337635 & 113815143225 & 29136676665600S-6\\
300 & 607500 & 13500 & 303600 & 303900 & 92172960000 & 23596277760000S-6\\
324 & 562500 & 13500 & 281088 & 281412 & 79010463744 & 20226678718464S-6\\
360 & 506250 & 13500 & 252945 & 253305 & 63981173025 & 16379180294400S-6\\
450 & 405000 & 13500 & 202275 & 202725 & 40915175625 & 10474284960000S-6\\
486 & 375000 & 13500 & 187257 & 187743 & 35065184049 & 8976687116544S-6\\
500 & 364500 & 13500 & 182000 & 182500 & 33124000000 & 8479744000000S-6\\
540 & 337500 & 13500 & 168480 & 169020 & 28385510400 & 7266690662400S-6\\
600 & 303750 & 13500 & 151575 & 152175 & 22974980625 & 5881595040000S-6\\
648 & 281250 & 13500 & 140301 & 140949 & 19684370601 & 5039198873856S-6\\
750 & 243000 & 13500 & 121125 & 121875 & 14671265625 & 3755844000000S-6\\
810 & 225000 & 13500 & 112095 & 112905 & 12565289025 & 3216713990400S-6\\
900 & 202500 & 13500 & 100800 & 101700 & 10160640000 & 2601123840000S-6\\
972 & 187500 & 13500 & 93264 & 94236 & 8698173696 & 2226732466176S-6\\
1000 & 182250 & 13500 & 90625 & 91625 & 8212890625 & 2102500000000S-6\\
1080 & 168750 & 13500 & 83835 & 84915 & 7028307225 & 1799246649600S-6\\
1250 & 145800 & 13500 & 72275 & 73525 & 5223675625 & 1337260960000S-6\\
1350 & 135000 & 13500 & 66825 & 68175 & 4465580625 & 1143188640000S-6\\
1458 & 125000 & 13500 & 61771 & 63229 & 3815656441 & 976808048896S-6\\
1500 & 121500 & 13500 & 60000 & 61500 & 3600000000 & 921600000000S-6\\
\hline
\end{tabular}
\end{center}
\begin{center}
\begin{tabular}{ccccccc}
\hline
$x$ & $y$ & $b=12\cdot 1125$ & $a$ & $d$ & $b^{2}$ & $b^{2}/a^{2}$ \\
\hline
1620 & 112500 & 13500 & 55440 & 57060 & 3073593600 & 786839961600S-6\\
1800 & 101250 & 13500 & 49725 & 51525 & 2472575625 & 632979360000S-6\\
1944 & 93750 & 13500 & 45903 & 47847 & 2107085409 & 539413864704S-6\\
2250 & 81000 & 13500 & 39375 & 41625 & 1550390625 & 396900000000S-6\\
2430 & 75000 & 13500 & 36285 & 38715 & 1316601225 & 337049913600S-6\\
2500 & 72900 & 13500 & 35200 & 37700 & 1239040000 & 317194240000S-6\\
2700 & 67500 & 13500 & 32400 & 35100 & 1049760000 & 268738560000S-6\\
2916 & 62500 & 13500 & 29792 & 32708 & 887563264 & 227216195584S-6\\
3000 & 60750 & 13500 & 28875 & 31875 & 833765625 & 213444000000S-6\\
3240 & 56250 & 13500 & 26505 & 29745 & 702515025 & 179843846400S-6\\
3750 & 48600 & 13500 & 22425 & 26175 & 502880625 & 128737440000S-6\\
4050 & 45000 & 13500 & 20475 & 24525 & 419225625 & 107321760000S-6\\
4500 & 40500 & 13500 & 18000 & 22500 & 324000000 & 82944000000S-6\\
4860 & 37500 & 13500 & 16320 & 21180 & 266342400 & 68183654400S-6\\
5000 & 36450 & 13500 & 15725 & 20725 & 247275625 & 63302560000S-6\\
5400 & 33750 & 13500 & 14175 & 19575 & 200930625 & 51438240000S-6\\
\hline
5832 & 31250 & 13500 & 12709 & 18541 & 161518681 & 41348782336S-6\\
6250 & 29160 & 13500 & 11455 & 17705 & 131217025 & 33591558400S-6\\
6750 & 27000 & 13500 & 10125 & 16875 & 102515625 & 26244000000S-6\\
7290 & 25000 & 13500 & 8855 & 16145 & 78411025 & 20073222400S-6\\
7500 & 24300 & 13500 & 8400 & 15900 & 70560000 & 18063360000S-6\\
\hline
8100 & 22500 & 13500 & 7200 & 15300 & 51840000 & 13271040000S-6\\
9000 & 20250 & 13500 & 5625 & 14625 & 31640625 & 8100000000S-6\\
9720 & 18750 & 13500 & 4515 & 14235 & 20385225 & 5218617600S-6\\
11250 & 16200 & 13500 & 2475 & 13725 & 6125625 & 1568160000S-6\\
12150 & 15000 & 13500 & 1425 & 13575 & 2030625 & 519840000S-6\\
12500 & 14580 & 13500 & 1040 & 13540 & 1081600 & 276889600S-6\\
\hline
\end{tabular}
\end{center}
%
\section{The Probability to get P322 Bounded Solutions.}
In this section our purpose is to calculate the probability to get solutions in the bounded interval for angle
$\theta = \arctan{a/b}$ bounded in the same P322 interval, $\theta \in (\pi/6, \pi/4)$, with respect to 
the total number of solutions for a given set of $Q$ values in interval $\theta \in (\pi/2, 0)$. 
\\\indent With base on the number of solutions for all $Q=1,2,\ldots,1125$:
$$
\frac{\text{Number of Solutions for } \theta \in (\pi/6, \pi/4)}
{\text{Number of Solutions for } \theta \in (\pi/2, 0)}
= \frac{3265}{50781} = 0.0642957,
$$
$$
\frac{\text{Number of Solutions for theta in P322 interval}}
{\text{Number of Solutions for } \theta \in (\pi/2, 0)} 
=\frac{3265-248}{50781} = \frac{3017}{50781} = 0.0594119848.
$$
In both cases the probability is 6\%.
\\\indent With base on the number of solutions for the same $Q$ used by the P322 scribe:
$$
\frac{\text{Number of Solutions for } \theta \in (\pi/6, \pi/4)}
{\text{Number of Solutions for } \theta \in (\pi/2, 0)}
=\frac{52}{614} = 0.08469,
$$
$$
\frac{\text{Number of Solutions for theta in P322 interval}}
{\text{Number of Solutions for } \theta \in (\pi/2, 0)} 
=\frac{51}{614} = 0.083.
$$
In both cases the probability is 8\%.
\\\indent With base on the number of different angles on the intervals:
\\\indent For all solutions:
\begin{equation*}
\frac{\text{Number of angles } \theta \in (\pi/6, \pi/4)}
{\text{Number of angles } \theta \in (\pi/2, 0)}
=\frac{382}{10378} = 0.0368,
\end{equation*}
$$
\frac{\text{Number of angles  theta in P322 interval}}
{\text{Number of angles } \theta \in (\pi/2, 0)}
=\frac{332}{10378} = 0.03199.
$$
Giving probabilities of 4\% and 3\%.
\\\indent For the same $Q$ values as in P322:
$$
\frac{\text{Number of angles } \theta \in (\pi/6, \pi/4)}
{\text{Number of angles } \theta \in (\pi/2, 0)}
=\frac{16}{193} = 0.0829,
$$
$$
\frac{\text{Number of angles  theta in P322 interval}}
{\text{Number of angles } \theta \in (\pi/2, 0)}
=\frac{15}{193} = 0.07772.
$$
In both cases the probability is again 8\%.
\begin{figure}[h!]
\centering
\caption{P322 histogram for angle $\theta$.}
\includegraphics[width=17cm]{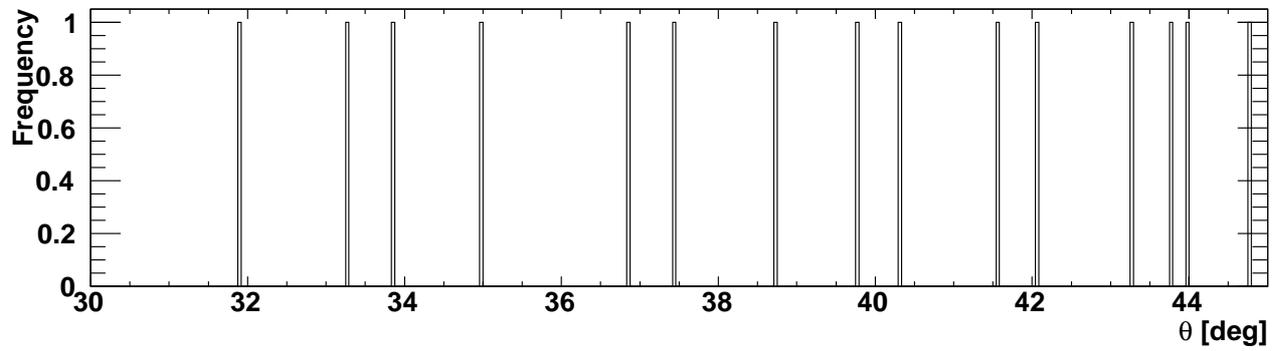}
\label{f:thetaP322}
\end{figure}
\begin{landscape}
\begin{figure}
\centering
\caption{Histograms obtained by using all the integers for $Q \in [1, 1125]$. Top is for $0 < \theta < \pi/2$.
Middle is for $\pi/3 < \theta < \pi/4$. Bottom is for $\pi/3 < \theta < \pi/4 \land \theta$ outside the P322 data interval.}
\includegraphics[width=21cm]{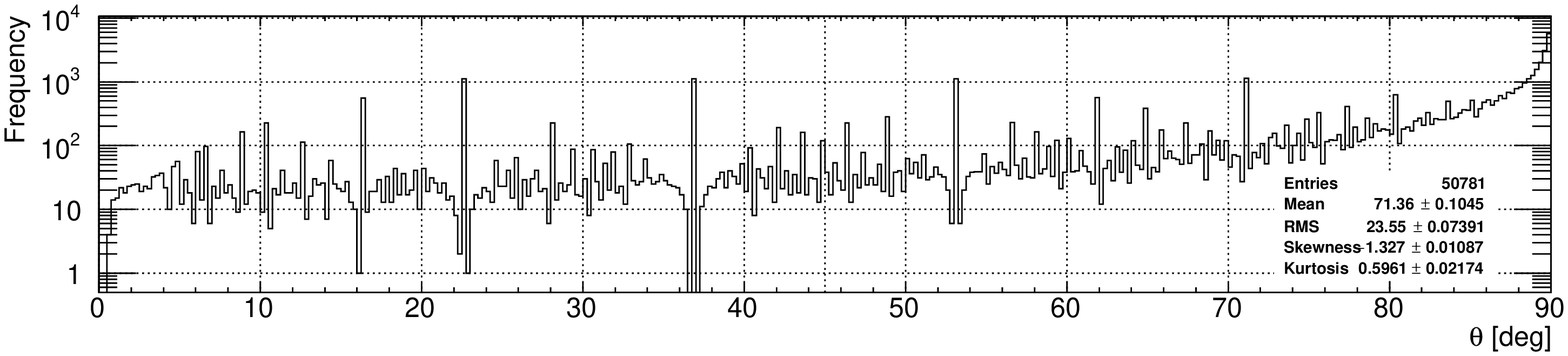}
\\\includegraphics[width=21cm]{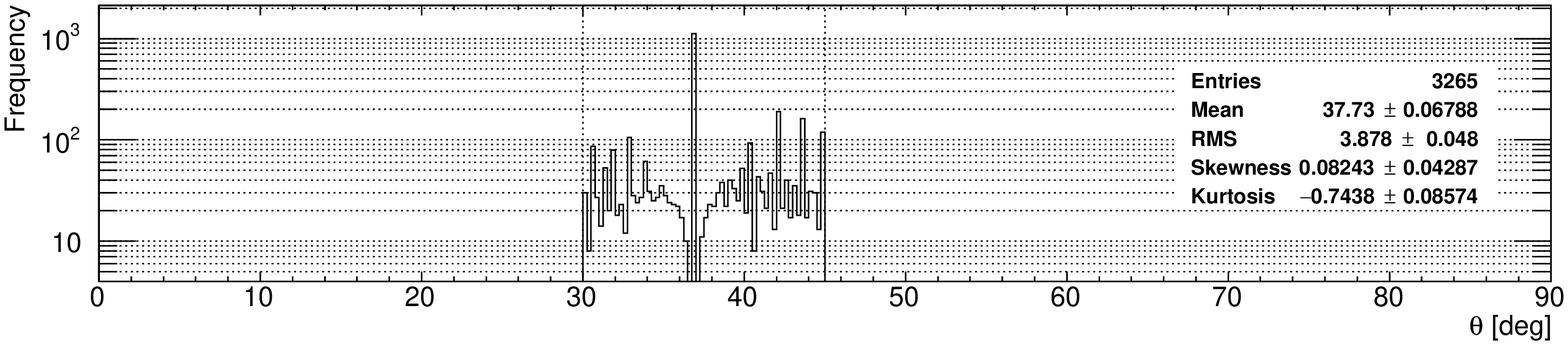}
\\\includegraphics[width=21cm]{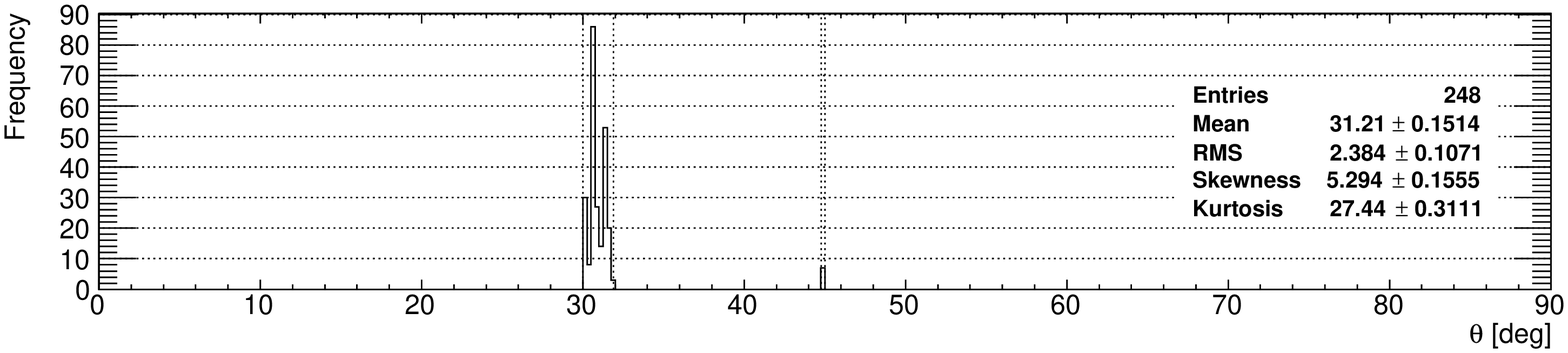}
\label{f:histoallQ}
\end{figure}
\begin{figure}
\centering
\caption{Histograms obtained by using the same values of $Q$ used by the scribe in P322. Top is for $0 < \theta < \pi/2$.
Middle is for $\pi/3 < \theta < \pi/4$. Bottom is for $\pi/3 < \theta < \pi/4 \land \theta$ outside the P322 data interval.}
\includegraphics[width=21cm]{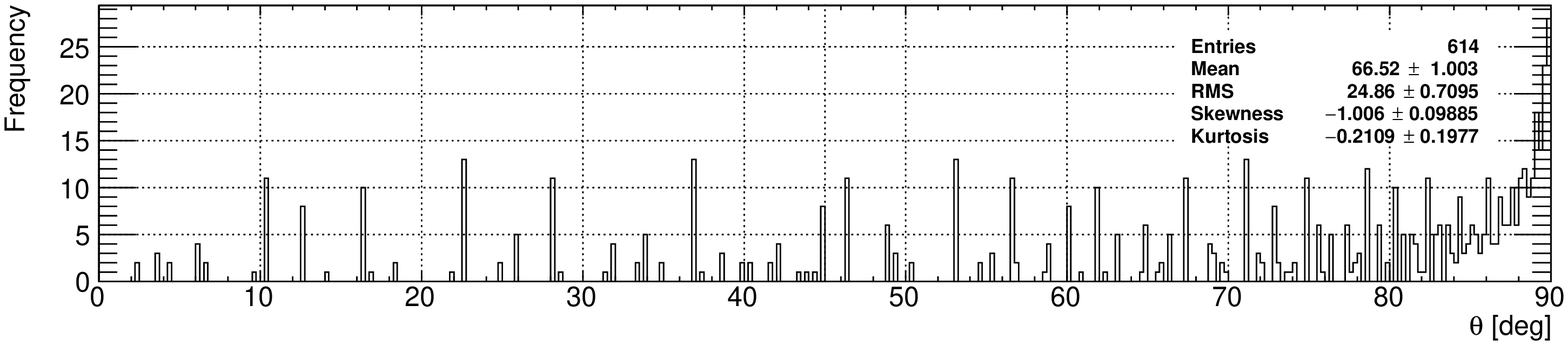}
\\\includegraphics[width=21cm]{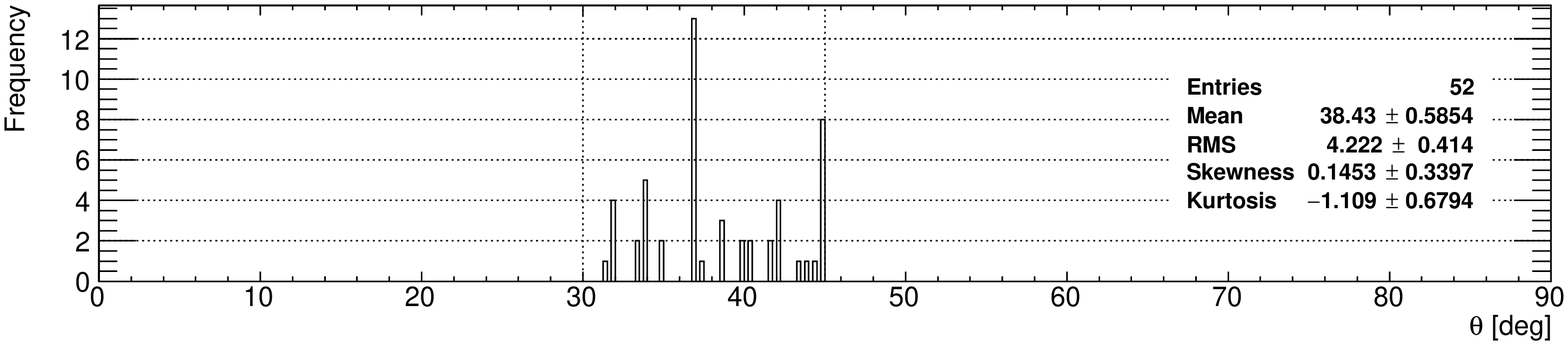}
\\\includegraphics[width=21cm]{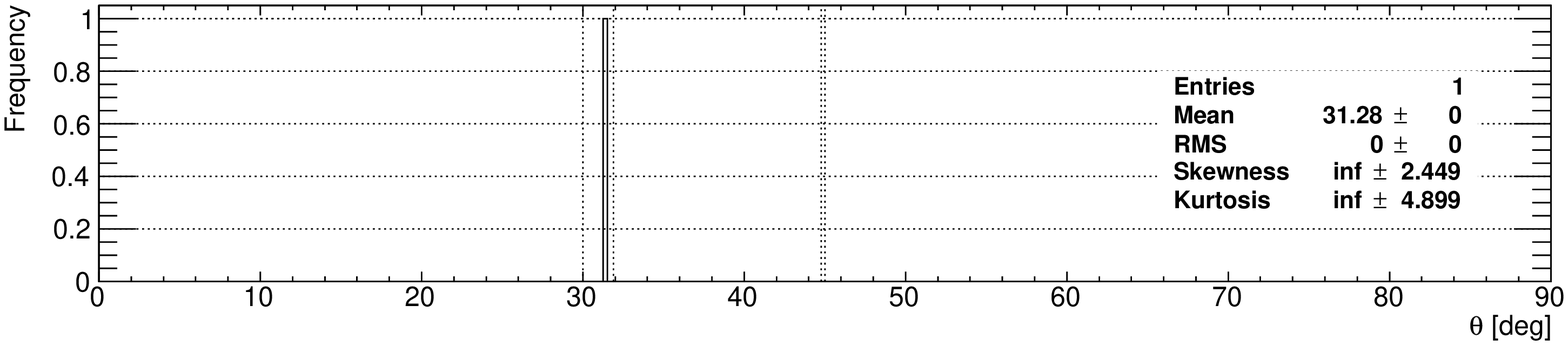}
\label{f:histoP322Q}
\end{figure}
\end{landscape}
\section{Discussion} 
\label{sec:Discussion}
We have presented a new algorithm to reconstruct the Plimpton 322 tablet, the algorithm is data driven and is straightforward,
also the factor twelve is justified inherent to the factorizing properties of the se\-xa\-ge\-si\-mal base.
Here we have applied the algorithm in three versions: 
\begin{enumerate}
\item when $Q=1$ by modifying the standard MT of reciprocals and using the factor $1.2$.
\item when we have generators $X$ and $Y$ with up to four sexagesits,
\item for integer solutions by using their scale generator $Q$ values since beginning, then we avoid to make a lot of operations 
to get $a$, $b$, $d$.
\end{enumerate}
The last version is the winner because use less operations than the other ones.
\subsection{P322 as a trigonometrical table.}
\begin{figure}
\centering
\caption{Representation of the carved information presented in tables YBC 7302 (left) and YBC 11120 (right). The number just above the internal circle $A$
is the circumference $c$, the number apart to the right side is $c^{2}$ and the number in the center is the area of the external circle $B$.}
\includegraphics[width=16cm]{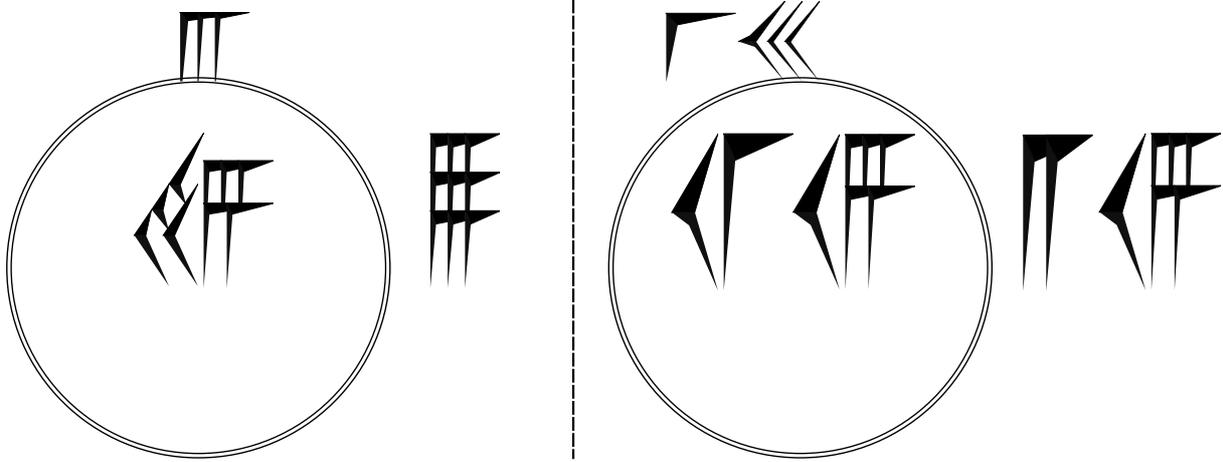}
\label{f:circles}
\end{figure}

Retractors of this theory, base their arguments in the modern point of view of a trigonometry based on the concept of a vector radius that sweep angles, 
and also that Old Babylonians did use the 
number $\pi$ as somehow equals to 3.
For the sake of discussion lets calculate the number $\pi$ in sexagesimal base.
\\\indent Up to $60^{-1}$
$$
\pi  \approx 03.~08 =_{10} 3.1(\overline{3}) = 3 + \frac{2}{15}.
$$
Up to $60^{-2}$ 
$$
\pi  \approx 03.~08~29 =_{10} 3.141(3\overline{8}) = 3 + \frac{509}{60^{2}}.
$$
\\Now lets perform the 
calculation of $\pi$ up to $60^{-3}$
\begin{equation}
\pi \approx 03.~08~29~44 =_{10} 3.141592(5926\ldots) = 3 + \frac{3823}{27000}
\label{eq:pi_3}
\end{equation}
this approximation gives $\pi$ with error $6.1\times 10^{-8}$ which is smaller than $60^{-3}$.
\\Up to $60^{-5}$ 
$$
\pi  \approx 03.~08~29~44~00~47 =_{10} 3 + \frac{28800 \cdot 3823 + 47}{60^{5}} = 3 + \frac{110102447}{777600000}.
$$
\\Up to $60^{-6}$ 
$$
\pi  \approx 03.~08~29~44~00~47~25 =_{10} 3 + \frac{60(28800 \cdot 3823 + 47)+25}{60^{6}} = 3 + \frac{1321229369}{9331200000}.
$$
\\Up to $60^{-7}$ 
$$
\pi  \approx 03.~08~29~44~00~47~25~53 =_{10} = 3 + \frac{396368810753}{2799360000000}.
$$
\\By using 8 sexagesits after the floating point:
\begin{equation}
\pi \approx 03.~08~29~44~00~47~25~53~07 =_{10} 3.14159265358979(0\ldots) = 3 + \frac{23782128645187}{167961600000000}
\label{eq:pi_8}
\end{equation}
gives $\pi$ with error $3\times 10^{-15}$ which is approx. $\frac{1}{2}\epsilon_{m}$.
\\\\\indent Thus we cannot underestimate the Old Babylonian power of calculation.
Performing $60^{6}$ in the CASIO fx-991 presents overflow, and in order to 
do calculations on computers with C++ without overflow, we need integers of at least 8 bytes of the type
\textsf{unsigned long long}.
\\\\\indent As a matter of fact the number $\pi$ appears in calculations of the circumference of the circle, and
its area. If the diameter of the circle is $d$, the circumference or perimeter is
$c = \pi d$, and its square is $c^{2} = \pi^{2} d^{2}$, then the area of the circle is
$A = \frac{1}{4} \pi d^{2} = \frac{1}{4} \, c d = \frac{1}{\pi}\frac{1}{4} c^{2}$,
but the information carved on tables corresponds to a from above approximate area that we have called $B$ 
\begin{equation}
B = \frac{1}{12} \, c^{2},
\label{e:areaB}
\end{equation}
and we have the real area is:
\begin{equation}
A = \frac{3}{\pi} \, B.
\label{e:areaA}
\end{equation}
In High Ener\-gy Physics we use Natural Units to describe the fundamental constants of nature but putting them equals to 1.
For example the speed of light and the Planck constant are both equal to 1.
This technique simplify the calculations and in the end the physicist convert its results to a desired system of units by multiplying them by 
nice factors.
In the same way if we assume the Old Babylonians perhaps use a similar technique by making the number 
$\pi$ equals to 3 and in the end the exact result is obtained by multiplyng by a nice approximation of $3/\pi$; 
but this would be denied because
not tables with this procedure have been found. 
But instead of that we can consider as in P322 the Old Babylonians avoid to deal with angles as multiples of $\pi$ and also to divide numbers in prime numbers
or multiples of prime numbers,
maybe because they hate to do approximations and they really like to do exact\footnote{With a finite number of sexagesits.} calculations.
Perhaps some day archeologists will find a table containing an approximation of $3/\pi = 00.~57~17~44~48~22\ldots$ and multiplications with it, 
for instance we can consider the use of the approximation
(\ref{e:areaB}) of (\ref{e:areaA}) expressed in the two circles carved on tables YBC 7302 and YBC 11120. 
See \citep{AMS-2} also Figure~\ref{f:circles}.
\\\\\indent We can note if the diameter of the circle $A$ with area $A$ (\ref{e:areaA}) and circumference $c$ is $d$,
then the circle $B$ with area $B$ (\ref{e:areaB}) will have a diameter 
\begin{equation}
d_{B} = \sqrt{\frac{\pi}{3}} \, d = (01.~01~23~58~34~08\ldots) \times d.
\label{e:extring}
\end{equation}
Thus the pictogram representation of concentric circles $A$ and $B$ will have 
the same concentric two circles pattern that we can see in tables YBC 7302 and YBC 11120.
Nevertheless from first principles the approximated area of a sector of a circle does not enter in calculations of angles, but
the circumference does it and the Old Babylonians take the value of the circumference and its squared in an exact way.
\\\\\indent Therefore in this work we consider Old Babylonians consciently avoid the use of $\pi$ in the calculation of areas and (if you allow) operations on angles
in the sense of the ratios of the squared sides of right triangles.
This is the reason that is so difficult to get prove the Old Babylonians did know the numerical value of the number $\pi$
in a better approximation than 3. On the other hand, accepting $\pi = 3$ for Old Babylonians results in a contradiction 
for a culture that did calculations with 8 (or 9) sexagesits of accuracy and does not give any explanation about the two concentric circles pattern found
in tablets YBC 7302 and YBC 11120.
\subsection{Prime Numbers.}
Once corrected the scribe errors indroduced in the table, and in special the 
prime $a=541$ in Row 9, we have the observation did by \citep[p. 340]{Buck1980} was:
\emph{``In the first 20000 integers there are about 2300 primes, which is about 10 percent; among 15 integers, selected at random from this interval, we
might, then, expect to see one or two primes, but certainly not eight!''}
\\\\\indent The explanation of this feature of P322 is the following: 
Considering the equivalence classes modulo 60 is well known that the primes greater than 5 will lie in the set
$C$ = \{ \cseven, \celeven, \cthirteen, \cseventeen, \cnineteen, \ctwentythree, \ctwentynine, \cthirtyone, 
\cthirtyseven, \cfourtyone, \cfourtythree, \cfourtynine, \cfiftythree, \cfiftynine, \cone \}.
In sexagesimal base the congruence classes modulo 60 is determined easily by the value of the last sexagesit.
Writing the diagonal theorem in terms of the congruence classes modulo $60$ in the form 
$[a]^{2} + [b]^{2} = [d]^{2}$ we got for the fifteen triangles:
\begin{eqnarray*}
\cfiftynine^{2} + \czero^{2} & = & \cfourtynine^{2} \\
\cseven^{2} + \cthirtysix^{2} & = & \ctwentyfive^{2} \\
\cfourtyone^{2} + \czero^{2} & = & \cfourtynine^{2} \\
\cfourtynine^{2} + \czero^{2} & = & \cone^{2} \\
\cfive^{2} + \ctwelve^{2} & = & \cthirtyseven^{2} \\
\cnineteen^{2} + \czero^{2} & = & \cone^{2} \\
\celeven^{2} + \czero^{2} & = & \cone^{2} \\
\cnineteen^{2} + \czero^{2} & = & \cfourtynine^{2} \\
\cone^{2} + \czero^{2} & = & \cfourtynine^{2} \\
\cfourtyone^{2} + \czero^{2} & = & \cone^{2} \\
\cfourtyfive^{2} + \czero^{2} & = & \cfifteen^{2} \\
\cfiftynine^{2} + \czero^{2} & = & \cfourtynine^{2} \\
\cfourtyone^{2} + \czero^{2} & = & \cfourtynine^{2} \\
\cthirtyone^{2} + \czero^{2} & = & \cfourtynine^{2} \\
\ctwentyeight^{2} + \cfourtyfive^{2} & = & \cfiftythree^{2} .
\end{eqnarray*}
Counting $[a]$ and $[d]$ occurrences of members of set $C$ we got 25 of 30, and the prime numbers came from
$d$ are according to Gauss of the form $4n+1$, prime numbers came from $a$ are of the form $4n+1$ and $4n+3$. 
We also note there are other congruences $\mod 60$ not in set $C$ that are in P322 as 
 \cfive, \cfourtyfive, \cfifteen, \ctwentyfive, \ctwentyeight.
Therefore P322 is not strictly a prime numbers generator algorithm, but certainly here there is an abundance of prime numbers and congruences elements of $C$
because of the chosen triples. 
For $[b]$ we have: \czero, \cthirtysix, \ctwelve, \cfourtyfive.
\\\\\indent In Section \ref{s:solutions} we observe several triangles with the same angles, and we can see the P322 scribe did prefer
to write the solutions in terms of the congruence classes modulo 60 instead of the triples that are multiple of these ones.
For examble when we use the factor $Q=06~40$ we can get the solutions for factors $03~20$, $01~20$, $50$, $40$, $20$, and  $5$; and the scribe
did prefer the simplified triples (Pythagorean).

\subsection{Pythagorean Triples}
Pythagorean triples are of the form $(a,b,d)$ with $a,b,d$ coprimes between them.
Threrefore when we simplify all the integer solutions for the 12-factor algorithm triples we noted several triangles are similar 
and thus they correspond to the same Pythagorean triple. We can see this in histograms of Figures~\ref{f:histoallQ} and \ref{f:histoP322Q}.
\\\\\indent Analyzing the possible called forgotten 16th row by the scribe, 
we have for the integer solutions scheme: 
\\$(40~00,\: 57~36,\: 67~24)$ with $Q=04~48$, $X=32~24$, $01~42~24$, this $Q$ was already used by the scribe;
for the scheme of generators $X$, $Y$, with up to 4 sexagesits: $(02~45,\: 04~48,\: 05~37)$ with $Q=24$, $X=02~42$, $Y=08~32$ and congruence
equation $\cfourtyfive^{2} + \cfourtyeight^{2} = \cthirtyseven^{2}$, this solution is also found in the 
integer solutions scheme;
for the $Q=1$ scheme $(07.~27~30,\: 12,\: 14.~12~30)$ with $X=06.~45$, $Y=21.~10$.
\\Here there are three options: we can simply say the scribe forgot to write this solution because he already used the same $Q=04~48$, 
or he discarded the solution $(40~00,\: 57~36,\: 67~24)$ because he wanted ``Pithagorean Triples'',
or forgot to calculate it because he didn't use the $Q=24$ scale factor.  
\\\\\indent The fiftenn row is saved from removal because the numbers of the triple in sexagesimal base and with the Old Babylonian place value notation 
looks like if $d$ were the prime $53$. 
\subsection{Bounded and Selected Solutions Criterion} We can get the 15 solutions presented in P322 by requiring to apply the integer solutions scheme 
for $\pi/6 < \theta < \pi/4$ with scale generators $Q$ in P322, see Table~\ref{t:QQQ}, such that $a$, $b$, $d$ be coprimes.
This additional cut reduce our prior probability error to get bounded solutions from 1/52 $\approx$ 2\%, 
see plots Figures~\ref{f:histoP322Q} and \ref{f:thetaP322}, eliminating 
the 16th triple $(40~00,\: 57~36,\: 67~24) =_{10} (2400, 3456, 4044)$ and now we get our new error to generate the P322 triples
to 0/51 = 0\%. 
\section{Conclusions}
The algorithm of factor 12 has unveiled the internal structure of the Old Babylonian sexagesimal ratio based trigonometry,
and the integer solutions scheme have been proved that is the most straightforward and plausible algorithm between all that claim to reconstruct
the P322 data in its own historical context and is the one with the minimum of arithmetical operations.
\\\\\indent It could be the simbol for \emph{ki.} in the first column, would mean in this context,
\emph{multiple of 12}?.
\\\\\indent Although the 4th column values were derived from algebra and the geometrical pro\-per\-ti\-es of different right triangles
by calculating ratios of the square of their sides,
the entire 4th column, this is, the numerical values together with their title,
can be considered as the first \emph{trigonometric identity} recorded in the history of mankind.
\\\\\indent The solutions carved on P322 are bounded for $\theta \in (\pi/6, \pi/4)$ and selected by choosing solutions triples with $a$, $b$, $d$ coprimes 
or solution ($a\times 60$, $b\times 60$, $d \times 60$) with $a$, $b$, $d$ coprimes.
\\\\\indent We have shown pictorially-numerically by absurd reductio that Old Babylonians knew the value of $\pi$ in a be\-tter approximation than 3.
They use 3 for fast approximations from above of circle area, therefore they have another factor 12 in
$
c^{2} = 12\, B.
$
But we have to be conservative and wait until more tablets with the two circles pattern could be found to consider this theory as proved.
\\\\\indent The Old Babylonian triples obtained by the algorithm of factor 12 in the integer solutions scheme are of the form:
\begin{equation}
a = \frac{1}{2} \, \left( \frac{12^{2} \, Q^{2}}{x} - x \right), \quad b = 12\,Q, \quad d = \frac{1}{2} \, \left( \frac{12^{2} \, Q^{2}}{x} + x \right), 
\end{equation}
with $\frac{12^{2} \, Q^{2}}{x} > x$ and $x$, $Q$, $\frac{12^{2} \, Q^{2}}{x}$ positive integers, thus
they compounds a different set of solutions that given by the
Euclid formulas:
$$
a = m^{2} - n^{2}, \quad b = 2 mn, \quad d = m^{2} + n^{2}, 
$$
with $m > n$ and $m$, $n$ positive integers.
However the number of solutions are less with respect to the obtained by the modified Euclid formulas
$$
a = k(m^{2} - n^{2}), \quad b = 2 kmn, \quad d = k(m^{2} + n^{2}), 
$$
with $m > n$ and $m$, $n$, $k$ positive integers.

Although the bundle factor $M=12$ does not generate all the Pythagorean triples, we can change it for another 
selections, e.g. 1, 2, 3, 4, 5, etc. to get all them with the general algorithm:
\begin{equation}
a = \frac{1}{2} \, \left( \frac{M^{2} \, Q_{M}^{2}}{x} - x \right), \quad b = M \, Q_{M}, \quad d = \frac{1}{2} \, \left( \frac{M^{2} \, Q_{M}^{2}}{x} + x \right),
\end{equation}
with $\frac{M^{2} \, Q_{M}^{2}}{x} > x$ and
$x$, $Q_{M}$, $M$, $\frac{M^{2} \, Q_{M}^{2}}{x}$ positive integers.

\section*{References}
\bibliographystyle{plainnat} 
%

\end{document}